\documentclass[11pt]{article}
\usepackage[utf8]{inputenc}
\usepackage{hyperref}

\usepackage{latexsym,epsfig}
\usepackage{amssymb,amsmath,amsfonts}
\usepackage{exscale} 
\usepackage{ifthen}
\usepackage{longtable}
\usepackage{subfigure}
\usepackage{caption}
\usepackage{cite}
\usepackage{lscape}

\begin{document}
\bibliographystyle{gabialpha}
\protect\pagenumbering{arabic}
\setcounter{page}{1}
 
\newcommand{\Zt}{\rm}

\newcommand{\ba}{\begin{array}}
\newcommand{\ea}{\end{array}}
\newcommand{\pot}{{\cal P}}
\newcommand{\curv}{\cal C}
\newcommand{\ddt} {\mbox{$\frac{\partial  }{\partial t}$}}
\newcommand{\hl}{\sf}
\newcommand{\hd}{\sf}

\renewcommand{\d}{\mathrm{d}}
\newcommand{\e}{\mathrm{e}}
\newcommand{\Ad}{\mbox{\rm Ad}}
\newcommand{\Adsm}{\mbox{{\rm \scriptsize Ad}}}
\newcommand{\ad}{\mbox{\rm ad}}
\newcommand{\adsm}{\mbox{{\rm \scriptsize ad}}}
\newcommand{\diag}{\mbox{\rm Diag}}
\newcommand{\sect}{\mbox{\rm sec}}
\newcommand{\id}{\mbox{\rm id}}
\newcommand{\idsm}{\mbox{{\rm \scriptsize id}}}
\newcommand{\eps}{\varepsilon}
\newcommand{\Summ}{P}

\newcommand{\aL}{\mathfrak{a}}
\newcommand{\bL}{\mathfrak{b}}
\newcommand{\mL}{\mathfrak{m}}
\newcommand{\kL}{\mathfrak{k}}
\newcommand{\gL}{\mathfrak{g}}
\newcommand{\nL}{\mathfrak{n}}
\newcommand{\hL}{\mathfrak{h}}
\newcommand{\pL}{\mathfrak{p}}
\newcommand{\uL}{\mathfrak{u}}
\newcommand{\lL}{\mathfrak{l}}

\newcommand{\kG}{{\tt k}}
\newcommand{\nG}{{\tt n}}

\newcommand{\Cart}{$G=K e^{\overline{\aL^+}} K$}
\newcommand{\Area}{\mbox{Area}}
\newcommand{\Hd}{\mbox{\rm Hd}}
\newcommand{\Hdim}{\mbox{\rm dim}_{\mbox{\rm \scriptsize Hd}}}
\newcommand{\Tr}{\mbox{\rm Tr}}
\newcommand{\bs}{{\cal \beta}}
\newcommand{\bv}{{\cal B}}

\newcommand{\nc}{{\cal N}}
\newcommand{\MM}{{\cal M}}
\newcommand{\Ch}{{\cal C}}
\newcommand{\clCh}{\overline{\cal C}}
\newcommand{\Sh}{\mbox{Sh}}
\newcommand{\smSh}{\mbox{{\rm \scriptsize Sh}}}
\newcommand{\Cnt}{\mbox{\rm C}}
\newcommand{\preim}{(\pi^F)^{-1}}

\newcommand{\NN}{\mathbb{N}} \newcommand{\ZZ}{\mathbb{Z}}
\newcommand{\QQ}{\mathbb{Q}} \newcommand{\RR}{\mathbb{R}}
\newcommand{\KK}{\mathbb{K}} \newcommand{\FF}{\mathbb{F}}
\newcommand{\CC}{\mathbb{C}} \newcommand{\EE}{\mathbb{E}}
\newcommand{\XX}{X}
\newcommand{\HH}{I\hspace{-2mm}H}
\newcommand{\norm}{\Vert\hspace{-0.35mm}|}
\newcommand{\Sph}{\mathbb{S}}
\newcommand{\ganz}{\overline{\XX}}
\newcommand{\rand}{\partial\XX}
\newcommand{\prodrand}{\partial\XX_1\times\partial\XX_2} 
\newcommand{\regrand}{\partial\XX^{reg}}
\newcommand{\singrand}{\partial\XX^{sing}}
\newcommand{\Frand}{\partial^F\XX}
\newcommand{\Lim}{L_\Gamma}          
\newcommand{\Flim}{F_\Gamma}
\newcommand{\reglim}{L_\Gamma^{reg}}
\newcommand{\radlim}{L_\Gamma^{rad}}
\newcommand{\raylim}{L_\Gamma^{ray}}
\newcommand{\horinf}{\mbox{Vis}^{\infty}}
\newcommand{\horF}{\mbox{Vis}^F}
\newcommand{\Sml}{\mbox{Small}}
\newcommand{\SmlF}{\mbox{Small}^F}

\newcommand{\ifl}{\qquad\Longleftrightarrow\qquad}
\newcommand{\at}{\!\cdot\!}
\newcommand{\ging}{\gamma\in\Gamma}
\newcommand{\xo}{{o}}
\newcommand{\gamo}{{\gamma\xo}}
\newcommand{\gam}{\gamma}
\newcommand{\gax}{h}
\newcommand{\gxi}{{G\!\cdot\!\xi}}
\newcommand{\bd}{$(b,\theta)$-densit}
\newcommand{\bt}{$(b,\theta)$-densit}
\newcommand{\cd}{$(\alpha,\Gamma\at\xi)$-density}
\newcommand{\be}{\begin{eqnarray*}}
\newcommand{\ee}{\end{eqnarray*}}

\newcommand{\an}{\ \mbox{and}\ }
\newcommand{\as}{\ \mbox{as}\ }
\newcommand{\diam}{\mbox{diam}}
\newcommand{\is}{\mbox{Is}}
\newcommand{\Ax}{\mbox{Ax}}
\newcommand{\Fix}{\mbox{Fix}}
\newcommand{\Par}{F}
\newcommand{\Min}{\mbox{Fix}}
\newcommand{\rel}{\mbox{Rel}_\Gamma}
\newcommand{\vol}{\mbox{vol}}
\newcommand{\Td}{\mbox{Td}}
\newcommand{\piF}{\pi^F}

\newcommand{\for}{\ \mbox{for}\ }
\newcommand{\pr}{\mbox{pr}}
\newcommand{\sh}{\mbox{sh}}
\newcommand{\shi}{\mbox{sh}^{\infty}}
\newcommand{\rank}{\mbox{rank}}
\newcommand{\supp}{\mbox{supp}}
\newcommand{\mass}{\mbox{mass}}
\newcommand{\kernel}{\mbox{kernel}}
\newcommand{\st}{\mbox{such}\ \mbox{that}\ }
\newcommand{\Stab}{\mbox{Stab}}
\newcommand{\Root}{\Sigma}
\newcommand{\Cone}{\mbox{C}}
\newcommand{\wrt}{\mbox{with}\ \mbox{respect}\ \mbox{to}\ }
\newcommand{\where}{\ \mbox{where}\ }

\newcommand{\thet}{\widehat H}

\newcommand{\con}{{\sc Consequence}\newline}
\newcommand{\rem}{{\sc Remark}\newline}
\newcommand{\prf}{{\sl Proof}}
\newcommand{\qed}{$\hfill\Box$}

\newenvironment{rmk} {\newline{\sc Remark.\ }}{}  
\newenvironment{rmke} {{\sc Remark.\ }}{}  
\newenvironment{rmks} {{\sc Remarks.\ }}{}  
\newenvironment{nt} {{\sc Notation}}{}  

\newtheorem{satz}{\bf Theorem}

\newtheorem{df}{\sc Definition}[section]
\newtheorem{cor}[df]{\sc Corollary}
\newtheorem{thr}[df]{\bf Theorem}
\newtheorem{lem}[df]{\sc Lemma}
\newtheorem{prp}[df]{\sc Proposition}
\newtheorem{ex}{\sc Example}
\newenvironment{pros}{{\sc Properties:}}


\title{{\sc Generalized conformal densities for higher products of rank one Hadamard spaces}}
\author{\sc Gabriele Link} 
\date{\today}
\maketitle
\begin{abstract} Let $\XX$ 
be a product of locally compact rank one Hadamard spaces and $\Gamma$ a  discrete group of isometries 
which contains two  elements projecting to a pair of independent rank one isometries in each factor.  
In \cite{1308.5584} we gave a precise description of the structure of the geometric limit set $\Lim$ of $\Gamma$; our aim in this paper is to describe this set from 
a measure theoretical point of view, using as a basic tool the properties of the exponent of growth of $\Gamma$ established in the aforementioned article.  We first show that the conformal density obtained from the classical Patterson-Sullivan construction is supported in a unique $\Gamma$-invariant subset of the geometric limit set; generalizing this classical construction 
we then obtain measures supported in each $\Gamma$-invariant subset of the regular limit set and investigate their properties. 

We remark that apart from 
Kac-Moody groups over finite fields acting on the Davis complex of their associated twin building, the probably most interesting examples to which our results apply are isometry groups of reducible CAT$(0)$-cube complexes without Euclidean factors. 

\end{abstract}

\vspace{0.2cm}

\section{Introduction}

Let $(\XX,d)$ be a product of $r$ locally compact Hadamard spaces $(\XX_i,d_i)$ endowed with the $\ell^2$-metric, which makes $\XX$ itself 
a locally compact Hadamard space, i.e. a locally compact complete simply connected metric spaces of non-positive Alexandrov curvature. 
It is well-known that every locally compact Hadamard space can be compactified 
by adding its  geometric boundary $\rand$ endowed with the cone topology (see \cite[Chapter II]{MR1377265}). If $\XX$ is a product space, then the {\hl regular geometric boundary} $\regrand$ of $\XX$ -- which consists of the set of equivalence classes of geodesic rays which do not project to a point in one of the factors -- is a 
dense open subset of $\rand$ homeomorphic to the Cartesian product of the geometric boundaries $\rand_i$ of the factors $\XX_i$ (which we call the Furstenberg boundary $\Frand$ of $\XX$) times a factor $E^+=\{\theta\in\RR^r_{>0}:\Vert \theta\Vert=1\}$; the projection to the last factor will be called the {\hl slope} of a 
point in $\regrand$. 
Every point $\tilde\eta$ in the {\hl singular geometric boundary}  $\singrand=\rand\setminus\regrand$ similarly has a well-defined slope $\theta=(\theta_1,\theta_2,\ldots,\theta_r)\in E:=\{\theta\in \RR_{\ge 0}^r: \Vert\theta\Vert=1\}$ \st $\theta_i=0$ for at least one $i\in\{1,2,\ldots,r\}$; notice that in this case the projection of $\tilde\eta$ to $\rand_i$ is not well-defined and $\tilde\eta$ is already completely determined by its slope and its projections $\eta_i$ to the geometric boundaries $\rand_i$ of the factors $\XX_i$ for which $\theta_i>0$.

For a group $\Gamma\subset \is(\XX_1)\times \is(\XX_2)\times\cdots\times\is(\XX_r)$ acting properly discontinuously by isometries on $\XX$ the limit set is defined by 
$\Lim:=\overline{\Gamma\at x}\cap\rand$, where $x\in\XX$ is arbitrary. In order to relate the critical exponent of a Fuchsian group to the Hausdorff 
dimension of its limit set, S.~J.~Patterson (\cite{MR0450547})  and D.~Sullivan (\cite{MR556586}) developed a theory of conformal densities. It 
turned out that for higher rank symmetric spaces and Euclidean buildings these densities in general detect only a small part of the geometric limit 
set (see \cite{MR1675889}). In order to measure the limit set in each invariant subset of the limit set, a class of generalized conformal densities 
were independently introduced in  \cite{MR1935549} and \cite{MR2062761}.  One of the main goals in this paper is to adapt this construction to 
 discrete groups  $\Gamma\subset\is(\XX_1)\times\is(\XX_2)\times\cdots\times\is(\XX_r)$ which contain a pair of isometries projecting to independent rank one elements in each factor. 
Related questions were considered by M.~Burger (\cite{MR1230298}) for graphs of convex cocompact groups in a product of rank one symmetric spaces, and  
by F.~Dal'bo and I.~Kim (\cite{MR2478811}) for discrete isometry groups of a product of two Hadamard manifolds of pinched negative curvature.  %

One important class of examples satisfying our conditions are Kac-Moody groups  $\Gamma$ over a finite field which act by isometries on a product 
$\XX=\XX_1\times\XX_2$, the CAT$(0)$-realization of the associated twin building ${\cal B}_+\times{\cal B}_-$.  Indeed, there  exists an element 
$h=(h_1,h_2)$ projecting to a rank one element in each factor by  Remark 5.4 and the proof of Corollary~1.3 in \cite{MR2585575}. Moreover, the 
action of the Weyl group produces many  axial isometries $g=(g_1,g_2)$ with $g_i$ rank one and independent from $h_i$ for $i=1,2$. Notice that if 
the order of the ground field is sufficiently large, then $\Gamma\subset\is(\XX_1)\times\is(\XX_2)$ is an irreducible lattice (see e.g. 
\cite{MR1715140} and \cite{MR2485882}). 

Moreover, according to the Rank Rigidity Theorem (\cite[Theorem A]{MR2827012}) every finite-dimensional CAT$(0)$-cube complex $\XX$ admitting a group 
$\Gamma$ of automorphisms without fixed point in the geometric compactification of $\XX$ and without a rank one isometry possesses a 
convex $\Gamma$-invariant subcomplex which is a product of two unbounded cube subcomplexes; so one inductively gets a convex $\Gamma$-invariant subcomplex of $\XX$ which can be decomposed into a finite product of rank one Hadamard spaces. In particular, our results apply to reducible finite-dimensional CAT$(0)$-cube complexes without Euclidean factor and discrete isometry groups as above. 

Apart from these examples possible factors of $\XX$ include 
locally compact Hadamard spaces of strictly negative Alexandrov 
curvature (compare \cite{MR2478811} in the manifold setting). In this special case every non-elliptic and non-parabolic isometry in one of the factors 
is already a rank one element.  Prominent examples here which are already covered by the above mentioned results of J.~F.~ Quint  and the author are 
Hilbert modular groups acting as irreducible lattices on a product of hyperbolic planes and graphs of convex cocompact groups of rank one symmetric spaces (see also \cite{MR1230298}).

A central role throughout the paper is played by the exponent of growth of $\Gamma$ of given slope $\theta=(\theta_1,\theta_2,\ldots,\theta_r)\in E$ introduced in Section~7 of\cite{1308.5584}. To recall its definition we fix a point $x=(x_1,x_2,\ldots, x_r)$ in $\XX$,  $\eps>0$, $n\gg 1$  and consider the cardinality $N_\theta^\eps(n)$ of  the set
\begin{align*} \big\{ \gamma=(\gamma_1,\gamma_2,\ldots,\gamma_r)\in\Gamma: &\  0<d(x,\gamma x)<n,\\ 
&\  \Big| \frac{d_i(x_i,\gamma_i x_i)}{d(x,\gamma x)}-\theta_i\Big|<\eps\ \text{ for all }\  1\le i\le r\}\big\}.
\end{align*} 
This number counts all orbit points $\gamma x$ 
of distance less than $n$ to the point $x$ which in addition are {``}close" to a geodesic ray in the class of a boundary point with slope $\theta$.
\begin{df}\label{expgrowth}
The {\hd exponent of  growth} of $\,\Gamma$ of slope $\theta\in E$ is defined by
$$\delta_\theta(\Gamma):=\lim_{\eps\to 0}\left( \limsup_{n\to\infty}\frac{\ln N_\theta^\eps(n)}{n}\right).$$
\end{df}
The quantity $\delta_\theta(\Gamma)$ can be thought of as a function of $\theta\in E$ which describes the exponential growth rate of orbit points 
converging to limit points of slope $\theta$. It is an invariant of $\Gamma$ which carries more information than the critical exponent $\delta(\Gamma)$;
from Theorem~7.6 in \cite{1308.5584} (compare also \cite[Theorem~7.4]{MR2629900} in the case of only two factors) it follows that there exists a unique slope $\theta^*\in E$ \st the exponent of growth 
of $\Gamma$ is maximal for this slope and equal to the critical exponent $\delta(\Gamma)$.  

Our first result concerns the measures on the geometric boundary obtained by the classical Patterson-Sullivan construction. Analogous to the case  of
symmetric spaces or Euclidean buildings of higher rank we have the following result:\\[2mm]
 {\bf Theorem A}$\quad$ {\sl The Patterson-Sullivan construction produces a conformal density with support in a single $\,\Gamma$-invariant subset of the geometric 
limit set. Every point in its support has slope $\theta^*$ as above.} \\[2mm]
Thus in order to measure the remaining $\Gamma$-invariant subsets of the limit set, we need a more sophisticated construction. Inspired by the paper \cite{MR1230298} of M.~Burger we will consider densities with more degrees of freedom than the classical conformal density.   

Before we can state the remaining results we need more definitions. We fix a base point $\xo=(\xo_1,\xo_2,\ldots,\xo_r)\in\XX$. For $\theta\in E$ 
we  denote $\rand_\theta$ the set of points in the geometric boundary of slope $\theta$ and $I^+(\theta)=\{ i\in\{1,2,\ldots,r\}: \theta_i>0\}$. Then according to the remark at the end of the first paragraph the strata  $\rand_\theta$ is homeomorphic to the Cartesian product of the geometric boundaries $\rand_i$ with $i\in I^+(\theta)$; for $\theta\in E^+$ this is obviously the whole Furstenberg boundary $\Frand=\rand_1\times\rand_2\times\cdots\times\rand_r$. For a point $\tilde\eta\in\rand_\theta$ and $i\in I^+(\theta)$ we will denote $\eta_i\in\rand_i$ the projection to the factor $\rand_i$.
Moreover, if $i\in\{1,2,\ldots,r\}$ and $\eta_i\in\rand_i$ we let 
$\bs_{\eta_i}(\cdot, \xo_i)$ denote the Busemann function centered at $\eta_i$ based at $\xo_i$.

\begin{df}\label{bdensi}
Let $\MM^+(\rand)$ denote the cone of positive finite Borel measures
on $\rand$, $\theta=(\theta_1,\theta_2,\ldots,\theta_r)\in E$ and $b=(b_1,b_2,\ldots,b_r)\in \RR^r$ \st $b_i=0$ for all $i\in\{1,2,\ldots,r\}\setminus I^+(\theta)$. 
A $\,\Gamma$-invariant {\hl $(b,\theta)$-density} is a
map 
$$\begin{array}{ccc}\mu:\XX &\to &\MM^+(\rand)\\
x&\mapsto&\mu_x\end{array}$$ 
\st for any $x=(x_1,x_2,\ldots,x_r)\in\XX$ the following three properties hold:
\begin{flushleft}
\begin{tabular}{rl}
 {\rm (i)}\ &   $\emptyset\neq\supp(\mu_{x})\subset\Lim\cap \rand_\theta$,\\[1mm]
{\rm (ii)}\  & $\forall \,\gamma\in\Gamma\qquad \gamma_*\mu_x=\mu_{\gamma  x}$\footnotemark[1],\\[1mm]
{\rm (iii)}\ & $\forall\, \tilde\eta
\in\supp(\mu_\xo)$\\[2mm]
&
\hspace*{2cm} $ \qquad\displaystyle \frac{d\mu_x}{d\mu_\xo}(\tilde\eta)=\e^{b_1\bs_{\eta_1}(\xo_1,x_1)+b_2
\bs_{\eta_2}(\xo_2,x_2)+\cdots+b_r \bs_{\eta_r}(\xo_r,x_r)}.$
\end{tabular}
\footnotetext[1]{Here $\gamma_*\mu_x$ denotes the measure defined by $\gamma_*\mu_x(E)=\mu_x(\gamma^{-1}E)$ for any Borel set $E\subset\rand$}
\end{flushleft}
\end{df} 

Notice that if $\theta_i=0$ for some $i\in\{1,2,\ldots,r\}$, then for $\tilde\eta\in\rand_\theta$ the projection $\eta_i\in\rand_i$ is not defined; however, 
 the condition $b_i=0$ ensures that the exponent in (iii) is well-defined. Moreover, the conformal density from Theorem~A is a special case of such a density with support in $\rand_{\theta^*}$ and parameters $b_i=\delta(\Gamma)\cdot\theta_i^*$, $i\in\{1,2,\ldots,r\}$. 

We next give a criterion for the existence of a \bd y.\\[2mm]
{\bf Theorem~B}$\quad$ {\sl If $\theta\in E^+$ is such that $\delta_\theta(\Gamma)>0$, then there exists a \bd y for some parameters 
$b=(b_1,b_2,\ldots,b_r)\in\RR^r$.}\\[2mm]
Notice that according to Theorem~7.9 of \cite{1308.5584} we have $\delta_\theta(\Gamma)>0$ for $\theta$ in the relative interior of the intersection of the limit cone $\ell_\Gamma$ with the vector subspace of $\RR^r$ it spans. 
In Section~\ref{genPseries} we will give an explicit construction of the \bd y from Theorem~B above.


The following results  about $(b,\theta)$-densities in particular apply to any conformal density supported in a single $\Gamma$-invariant subset of the geometric limit set, not 
only the one obtained by the classical Patterson-Sullivan construction. Our 
main tool is a so-called shadow lemma for $(b,\theta)$-densities, which is a generalization of the well-known shadow lemma for conformal densities. It first gives a condition 
for the parameters of a $(b,\theta)$-density in terms of the exponent of growth. \\[2mm]
{\bf Theorem C}$\quad$ {\sl If a $\,\Gamma$-invariant $(b,\theta)$-density exists for some  $\theta=(\theta_1,\theta_2,\ldots,\theta_r)\in E^+$, then
\vspace{-1mm}
\[ \delta_\theta(\Gamma)\le \sum_{i=1}^r b_i\cdot\theta_i.\]}

The following subsets of the geometric limit set will play an important role in the sequel.
\begin{df}\label{raddef}
A point $\tilde \xi\in \rand_\theta$  is called a {\hd radial limit point} of $\,\Gamma$ if there exists a sequence $(\gamma_n)=\big((\gamma_{n,1},\gamma_{n,2},\ldots,\gamma_{n,r})\big)\subset\Gamma$ 
\st $\gamma_n\xo$ converges to $\tilde\xi$ and \st for all $i\in I^+(\theta)$ the sequence 
 $\gamma_{n,i}\xo_i\subset\XX_i$ stays at bounded distance of one (and hence any) 
geodesic ray in the class of $\xi_i\subset\rand_i$.

We will denote the set of all radial limit points of $\,\Gamma$ by $\radlim$. 
\end{df}
Notice that in general a radial limit point $\tilde \xi\in\rand$ is not  approached by a sequence $\gamma_n\xo\subset\XX$ 
staying at bounded distance of a geodesic ray in the class of $\tilde \xi$.   

Our next statement shows that for certain \bd ies the corresponding exponent of growth $\delta_\theta(\Gamma)$ is completely determined by the parameters $\theta\in E^+$ 
and $b=(b_1,b_2,\ldots,b_r)\in\RR^r$:\\[2mm] 
{\bf Theorem D}$\quad$ {\sl If $\theta=(\theta_1,\theta_2,\ldots,\theta_r)\in E^+$, and $\mu$ is a $\,\Gamma$-invariant  \bd y which gives positive measure to the radial limit set, then 
\[ \delta_\theta(\Gamma)=\sum_{i=1}^r b_i\cdot \theta_i.\]}
The following theorem gives a restriction for  the atomic part of our measures.\\[2mm] 
{\bf Theorem E}$\quad$ {\sl If $\theta\in E^+$ \st  $\delta_\theta(\Gamma)>0$, and $\mu$ is a $\,\Gamma$-invariant  \bd y, then a radial limit point is not a point mass for $\mu$. 
}\\[-1mm]

%
%
%

The paper is organized as follows: 
In Section~\ref{Prelim} we recall  basic facts about Hadamard spaces and rank one isometries. Section~\ref{prodHadspaces} deals with the product case and provides some tools for the proof 
of the so-called shadow lemma in Section~7.  In Section~\ref{ExpGrowth} we introduce and study the properties of the exponent of growth. Section~\ref{ClasPatSulconst} recalls the 
classical Patterson-Sullivan construction in our setting. The main new result here is  Theorem A. In Section~\ref{genPseries} we introduce a generalized Poincar{\'e} series that  
allows to construct \bd ies, and therefore proves Theorem B. Using the shadow lemma, in Section~\ref{Propbdies} we deduce properties of \bd ies and prove Theorems C, D and E. 
\\[1mm]

{\bf Acknowledgements:}  The first draft of this paper was initiated during the author's stay at IHES in Bures-sur-Yvette. She warmly thanks the institute for its hospitality and the inspiring atmosphere.

\section{Preliminaries}\label{Prelim}

The purpose of this section is to introduce terminology and notation and to summarize basic results about Hadamard spaces and rank one isometries. 
The main references here  are \cite{MR1744486} and  \cite{MR1377265} (see also \cite{MR1383216}, and \cite{MR823981},\cite{MR656659} in the case of Hadamard manifolds). 

Let $(\XX,d)$ be a metric space. A {\hl geodesic path} joining $x\in\XX$ to $y\in\XX$  is a map $\sigma$ from a closed interval $[0,l]\subset \RR$ to $\XX$ \st $\sigma(0)=x$, 
$\sigma(l)=y$ and $d(\sigma(t), \sigma(t'))=|t-t'|$ for all $t,t'\in [0,l]$.  We will denote such a geodesic path $\sigma_{x,y}$. $\XX$ is called {\hl geodesic} if any two points 
in $\XX$ can be connected by a geodesic path, if this path is unique we say that $\XX$ is {\hl uniquely geodesic}. In this text $\XX$ will be a Hadamard space, i.e. a complete 
geodesic metric space in which all triangles satisfy the CAT$(0)$-inequality. This implies in particular that $\XX$ is simply connected and uniquely geodesic.   A {\hl geodesic} 
or {\hl geodesic line} in $\XX$ is a map $\sigma:\RR\to\XX$ \st $d(\sigma(t), \sigma(t'))=|t-t'|$ for all $t,t'\in\RR$, a {\hl geodesic ray} is a map $\sigma:[0,\infty)\to \XX$ 
\st $d(\sigma(t), \sigma(t'))=|t-t'|$ for all $t,t'\in [0,\infty)$. Notice that in the non-Riemannian setting completeness of $\XX$ does not imply geodesically completeness, i.e. 
not every geodesic path or ray can be extended to a geodesic.


From here on we will assume that $\XX$ is a locally compact Hadamard space.  The geometric boundary $\rand$ of
$\XX$ is the set of equivalence classes of asymptotic geodesic rays endowed with the cone topology (see e.g. \cite[chapter~II]{MR1377265}). The action of the isometry group 
$\is(\XX)$ on $\XX$ naturally extends to an action by homeomorphisms on the geometric boundary. Moreover, since $\XX$ is locally compact, this boundary $\rand$ is compact 
and the space $\XX$ is a dense and open subset of the compact space $\ganz:=\XX\cup\rand$.  For $x\in\XX$ and $\xi\in\rand$ arbitrary there exists a  geodesic ray emanating 
from $x$ which belongs to the class of $\xi$. We will denote such a ray $\sigma_{x,\xi}$.

We say that two points $\xi$, $\eta\in\rand$ can be joined by a geodesic if there exists a geodesic $\sigma:\RR\to\XX$ \st $\sigma(-\infty)=\xi$ and 
$\sigma(\infty)=\eta$. It is well-known that if $\XX$ is  CAT$(-1)$, i.e. of negative Alexandrov curvature bounded above by $-1$, then every pair of
distinct points in the geometric boundary can be joined by a geodesic. This is not true in the general CAT$(0)$-case. 

Let $x, y\in \XX$, $\xi\in\rand$ and $\sigma$ a geodesic ray in the class of $\xi$. We set 
\begin{equation}\label{buseman}
 \bs_{\xi}(x, y)\,:= \lim_{s\to\infty}\big(d(x,\sigma(s))-d(y,\sigma(s))\big).
\end{equation}
This number is independent of the chosen ray $\sigma$, and the function
\begin{align*} \bs_{\xi}(\cdot , y): \quad \XX &\to  \RR\\
x &\mapsto  \bs_{\xi}(x, y)\end{align*}
is called the {\hl Busemann function} centered at $\xi$ based at $y$ (see also \cite{MR1377265}, chapter~II). From the definition one immediately gets the following 
properties of the Busemann function:
\begin{align*}
|\bs_{\xi}(x, y)|&\le d(x,y)\\
\tag{anti-symmetry}\bs_{\xi}(x,y)&= -\bs_\xi(y,x)\\
\tag{cocycle identity} \bs_{\xi}(x, z)&=\bs_{\xi}(x, y)+\bs_{\xi}(y,z)\\
\tag{ Is$(\XX)$-invariance} \bs_{\xi}(x, y) &= \bs_{g\cdot\xi}(g\at x,g\at y) 
\end{align*}
for all $x,y,z\in\XX$, $\xi\in\rand$ and $g\in\is(\XX)$. Moreover, $\bs_{\xi}(x,y)=d(x,y)$ if and only if $y$ is a point on the geodesic ray $\sigma_{x,\xi}$, and we 
have the following easy 
\begin{lem}\label{distbusbounded}
Let $c>0$, $x,z\in \XX$ 
and $\xi\in\rand$ \st $d(z,\sigma_{x,\xi})<c$. Then 
$$0\le d(x,z)-\bs_\xi(x,z)<2c.$$
\end{lem}
\prf. \ \  The first inequality follows from $|\bs_{\xi}(x, y)|\le d(x,y)$. 
For the second one let $y\in\XX$ be a point on the geodesic ray $\sigma_{x,\xi}$ \st $d(z,y)<c$. Then for all $s>d(x,y)$ 
we have by the triangle inequality
\begin{align*} 
d(x,\sigma_{x,\xi}(s))-d(z,\sigma_{x,\xi}(s))& \ge  d(x,\sigma_{x,\xi}(s))-d(z,y)-d(y,\sigma_{x,\xi}(s))\\
&= d(x,y)-d(z,y)> d(x,y)-c,
\end{align*}
hence $\qquad d(x,z)-\bs_\xi(x,z)\le d(x,y)+c -d(x,y)+ c=2c$.\qed\\[-1mm]

A geodesic $\sigma: \RR\to\XX$  is said to {\hl bound a flat half-plane}  if there exists a closed convex subset $\iota([0,\infty)\times \RR)$ in $\XX$ isometric to 
$[0, \infty)\times\RR$ \st $\sigma(t)=\iota(0,t)$ for all $t\in \RR$. Similarly, a geodesic $\sigma: \RR\to\XX$  bounds a {\hl flat strip} of width $c>0$   if there 
exists a closed convex subset $\iota([0,c]\times \RR)$ in $\XX$ isometric to $[0, c]\times\RR$ \st $\sigma(t)=\iota(0,t)$ for all $t\in \RR$.
We call a geodesic $\sigma:\RR\to\XX$ a {\hl rank one geodesic} if  $\sigma$ does not bound a
flat half-plane. 

The following important lemma states that even though we cannot join any two distinct points in the geometric boundary of $\XX$, given a rank one geodesic we can at least join points in a neighborhood of its extremities. More precisely, we have the following well-known
\begin{lem}\label{joinrankone} (\cite{MR1377265}, Lemma III.3.1)\ 
Let $\sigma:\RR\to\XX$ be a rank one geodesic. Then there exist $c>0$ and 
neighborhoods $U^-$, $U^+ $ of $\sigma(-\infty)$, 
$\sigma(\infty)$ in $\ganz$ \st for any $\xi\in U^-$ and $\eta \in U^+$ there exists a rank one geodesic joining $\xi$ and $\eta$. For any such geodesic $\sigma'$ we have $d(\sigma', \sigma(0))\le c$.
\end{lem}

The following kind of isometries will
 play a central role in the sequel. 
\begin{df} \label{axialisos}
An isometry $h$ of $\XX$ is called {\hd axial}, if there exists a constant\break
$l=l(h)>0$ and a geodesic $\sigma$ \st
$h(\sigma(t))=\sigma(t+l)$ for all $t\in\RR$. We call
$l(h)$ the {\hd translation length} of $h$, and $\sigma$
an {\hd axis} of $h$. The boundary point
$h^+:=\sigma(\infty)$ is called the {\hd attractive fixed
point}, and $h^-:=\sigma(-\infty)$ the {\hd repulsive fixed
point} of $h$. We further set
$\Ax(h):=\{ x\in\XX:  d(x,h x)=l(h)\}$.
\end{df}

We remark that $\Ax(h)$ consists of the union of parallel geodesics
translated by $h$, and 
$\overline{\Ax(h)}\cap\rand$ is exactly the set of fixed points of
$h$.
\begin{df}
An axial isometry  is called {\hd rank one} if it possesses a rank one axis. Two rank one isometries are called 
{\hd independent}, if their fixed point sets are disjoint. 
\end{df}

Notice that if $h$ is rank one, then $h^+$ and $h^-$ are the only fixed points of $h$. 
Let us recall the north-south dynamics of rank one isometries.
\begin{lem}\label{dynrankone}(\cite{MR1377265}, Lemma III.3.3)\ 
Let $h$ be a rank one isometry. Then
\begin{enumerate}
\item[(a)] any $\xi\in\rand\setminus\{h^+\}$ can be joined to $h^+$ by a geodesic, and every geodesic joining  $\xi$ to $h^+$ is rank one, 
\item[(b)] given neighborhoods $U^-$ of $h^-$ and $U^+$ of $h^+$ in $\ganz$ 
there exists $N_0\in\NN$ \st\\
 $h^{-n}(\ganz\setminus U^+)\subset U^-$ and
$h^{n}(\ganz\setminus U^-)\subset U^+$ for all $n\ge N_0$.
\end{enumerate}
\end{lem}

If  $\Gamma$ is a group acting 
by isometries on a locally compact Hadamard space $\XX$ we define its {\hl geometric limit set}  by
$\Lim:=\overline{\Gamma\at x}\cap\rand$, where $x\in \XX$ is arbitrary.

\section{Products of Hadamard spaces}\label{prodHadspaces}

Now consider $r$ locally compact Hadamard spaces  $(\XX_1,d_1)$, $(\XX_2,d_2),\ldots, (\XX_r,d_r)$  and their Cartesian product $\XX=\XX_1\times \XX_2\times\cdots\times\XX_r$   endowed with the 
distance\break $d=\sqrt{d_1^2+d_2^2+\cdots+d_r^2}$. Notice that $(\XX,d)$ 
is again a locally compact Hadamard space.  

We denote $\RR_{\ge 0}^r:=\big\{(t_1,t_2,\ldots,t_r)\in\RR^r: t_i\ge 0\ \text{for all}\ i\in\{1,2,\ldots,r\}\big\}$ and\break
$\RR_{>0}^r:=\big\{(t_1,t_2,\ldots,t_r)\in\RR^r: t_i > 0\ \text{for all}\ i\in\{1,2,\ldots,r\}\big\}$. 
To any pair of points $x=(x_1,x_2,\ldots, x_r)$, $z=(z_1,z_2,\ldots, z_r)\in \XX$ we associate the vector 
\begin{equation}\label{distancevector}
H(x,z):= \left(\begin{array}{c} d_1(x_1,z_1)\\ d_2(x_2,z_2)\\\vdots\\ d_r(x_r,z_r)\end{array}\right)\in \RR_{\ge 0}^r\,,
\end{equation} 
which we call the {\hl distance vector} of the pair $(x,z)$. Notice that  if $\Vert \cdot\Vert$ denotes the Euclidean norm in $\RR^r$, we clearly have $\Vert H(x,z)\Vert =d(x,z)$. For $z\neq x$ we further define the {\hl direction} of $z$  \wrt $x$ by the unit vector
\begin{equation}\label{slopeangle}
\thet(x,z):=\frac{H(x,z)}{d(x,z)}  \in\RR_{\ge 0}^{r}\,.
\end{equation}


The following lemma is immediate and states that distance
 vectors and directions are invariant by $\is(\XX_1)\times\is(\XX_2)\times\cdots\times\is(\XX_r)$.
\begin{lem}\label{Ginvariance}
If $g=(g_1,g_2,\ldots,g_r)\in \is(\XX_1)\times\is(\XX_2)\times\cdots\times\is(\XX_r)$, $x=(x_1,x_2,\ldots,x_r)$, $z=(z_1,z_2,\ldots,z_r)\in\XX$, then 
$$ H(gx,gz)=H(x,z)\quad\mbox{ and }\quad \thet(gx,gz)=\thet(x,z).$$
\end{lem}

Denote $p_i:\XX\to \XX_i$, $i\in\{1,2,\ldots,r\}$, the natural projections. Every geodesic path $\sigma:[0,l]\to\XX\,$ can be written as a product $\sigma(t)=(\sigma_1(t\cdot \theta_1), \sigma_2(t\cdot \theta_2),\ldots, \sigma_r(t \cdot \theta_r))$, where  $\sigma_i$ are geodesic paths in $\XX_i$, $i=1,2,\ldots r$, and the $\theta_i\ge 0$ satisfy 
\[\sum_{i=1}^r \theta_i^2=1.\]  
The unit vector 
\[ \text{sl}(\sigma):= \left(\begin{array}{c} \theta_1\\ \theta_2\\\vdots\\ \theta_r\end{array}\right)\in E:=\{\theta\in  \RR_{\ge 0}^r :\ \Vert \theta\Vert =1\}\]
equals the direction of the points $\sigma(t)$, $t\in (0,l]$, with respect to $\sigma(0)$ and  is called the {\hl slope of $\sigma$}.  We say that a geodesic path $\sigma$ is  {\hl regular} if its slope does not possess a coordinate zero, i.e. if 
\[ \text{sl}(\sigma)\in E^+:=\{\theta\in  \RR_{> 0}^r :\ \Vert \theta\Vert =1\};\]
otherwise $\sigma$ is said to be {\hl singular}. In other words, $\sigma$ is regular if none of the projections $p_i(\sigma([0,l]))$, $i\in\{1,2,\ldots, r\}$,   is a point.

If  $x\in\XX$  and $\sigma:[0,\infty)\to\XX$ is an arbitrary geodesic ray, then by elementary geometric estimates one has the relation 
\[
 \text{sl}(\sigma)= \lim_{t\to\infty} \thet(x,\sigma(t)) =\lim_{t\to\infty} \frac{H(x,\sigma(t))}{d(x,\sigma(t))}
  \]
 between the slope of $\sigma$ and the directions of $\sigma(t)$, $t>0$, \wrt $x$. 
 Similarly, one can easily show that any two geodesic rays representing the same (possibly singular) point in the geometric boundary necessarily have the same slope. 
So  we may define the slope $\text{sl}(\tilde\xi)$ of a point $\tilde\xi\in\rand$ as the slope of an arbitrary geodesic ray representing $\tilde\xi$. The {\hl regular geometric boundary} $\regrand$ and  the {\hl singular geometric boundary} $\singrand$ of $\XX$ are then naturally defined by 
 \[ \regrand:=\{\tilde\xi\in\rand:\ \text{sl}(\tilde\xi)\in E^+\},\quad \singrand:=\rand\setminus\singrand;\]
the singular boundary $\singrand$ consists of equivalence classes of geodesic rays in $\XX$ which project to a point in at least one of the factors $\XX_i$.
More precisely, given\break $\theta=(\theta_1,\theta_2,\ldots,\theta_r)\in E$, we can  define the subset
\begin{equation}\label{randtheta}
 \rand_\theta:=\{\tilde\xi\in\rand:\ \text{sl}(\tilde\xi)=\theta\}
 \end{equation}
of the geometric boundary; we further denote $I^+(\theta):=\{i\in\{1,2,\ldots,r\}: \theta_i>0\}$. It is easy to see that two geodesic rays  $\sigma$, $\sigma'$  in $\XX$ represent the same point in $\rand_\theta$ if and only if  $\sigma_i(\infty)=\sigma_i'(\infty)$ for all $i\in I^+(\theta)$. Hence $\rand_\theta$ is homeomorphic to the Cartesian product of the geometric boundaries $\rand_i$ with $i\in I^+(\theta)$.

%


We further remark that a sequence $(y_n)=\big((y_{n,1},y_{n,2},\ldots,y_{n,r})\big)\subset\XX$ converges to a point $\tilde\eta\in\rand_\theta$ 
if and only if $y_{n,i}\to \eta_i$ for all $i\in I^+(\theta)$ and $\thet(x,y_n)\to\theta$ as $n\to\infty$ for some (and hence any) fixed $x\in\XX$. 

For higher rank symmetric spaces and Bruhat-Tits buildings there is a well-known notion of Furstenberg boundary which -- for a product of rank one spaces --
coincides with the product of the geometric boundaries. In our setting we choose to call the product $\rand_1\times \rand_2\times\cdots\times\rand_r$ endowed with
the product topology the {\hl Furstenberg boundary} $\Frand$ of $\XX$. Since $\regrand$ is homeomorphic to $\Frand\times E^+$ 
we have a natural projection
$$ \hspace{1cm} \ba{lccc}\pi^F\,: & \regrand&\to& \Frand\\
&(\xi_1,\xi_2,\ldots,\xi_r,\theta)&\mapsto & (\xi_1,\xi_2,\ldots,\xi_r)\, \ea $$
and a natural action of the  group $\is(\XX_1)\times\is(\XX_2)\times\cdots\times\is(\XX_r)$ by homeomorphisms on the Furstenberg boundary of $\XX$. 

We say that two points $\xi=(\xi_1,\xi_2,\ldots, \xi_r)$, $\eta=(\eta_1,\eta_2,\ldots, \eta_r)\in\Frand$ are {\hl opposite} if $\xi_i$ and $\eta_i$ can be joined by a geodesic in $\XX_i$ for all $i\in\{1,2,\ldots,r\}$. 

For $x=(x_1,x_2,\ldots,x_r)$ and $z=(z_1,z_2,\ldots,z_r)\in\XX$ 
\st $z_i\neq x_i$  for all\break $i\in\{1,2,\ldots,r\}$, the set 
\begin{align}\label{Weylchamberdef}
\Ch_{x,z}:=& \big\{(\sigma_{x_1,z_1}(t_1),\sigma_{x_2,z_2}(t_2),\ldots,\sigma_{x_r,z_r}(t_r))\in\XX : \nonumber\\
&\hspace*{4.4cm}  0\le t_i\le d_i(x_i,z_i)\quad\text{for}\quad i\in\{1,2,\ldots,r\}\big\}
\end{align}
is called the {\hl Weyl chamber} from $x$ to $z$. 
Notice that if $z_i=x_i$ for some $i\in\{1,2,\ldots,r\}$, then $\sigma_{x_i,z_i}$ is not defined, so the assignment in (\ref{Weylchamberdef}) 
is not well-defined. In this case we set $I(x,z):=\{i\in\{1,2,\ldots,r\}: x_i=z_i\}$ and define $\Ch_{x,z}$
\begin{align*}
\Ch_{x,z}:=& \big\{(y_1,y_2,\ldots,y_r)\in\XX :\ y_i\in\XX_i\quad\text{arbitrary if }\  i\in I(x,z),\\
&\hspace*{1.2cm}  y_i=\sigma_{x_i,z_i}(t_i)\quad\text{with}\quad 0 \le t_i\le d_i(x_i,z_i)\quad\text{if }\ i\in\{1,2,\ldots,r\}\setminus I(x,z)\big\}.
\end{align*}
We remark that in the degenerated case $z=x$ our definition gives $\Ch_{x,x}=\{x\}$.

Similarly, for $x=(x_1,x_2,\ldots,x_r)\in\XX$, $\theta\in E$ and $\tilde\xi\in\rand_\theta$
we call 
\begin{align}\label{Weylsingular}
\Ch_{x,\tilde\xi}:=& \big\{(y_1,y_2,\ldots,y_r)\in\XX :\ y_i=\sigma_{x_i,\xi_i}(t_i)\quad\text{with}\quad t_i\ge 0 \quad\text{if }\ i\in I^+(\theta), \nonumber\\
&\hspace*{3.9cm} y_i\in\XX_i\quad\text{arbitrary if }\  i\in \{1,2,\ldots,r\}\setminus I^+(\theta) \big\}
\end{align}
the Weyl chamber with apex $x$ in the class of $\tilde\xi$. 
In this way we have defined $\Ch_{x,z}$ for any $x\in\XX$ and $z\in\overline\XX $. Notice that while Weyl chambers in the class of a regular boundary point are homeomorphic to $\RR_{\ge 0}^r$, a Weyl chamber in the class of a singular boundary point 
is much bigger. 


The {\hl Weyl chamber shadow} of a set $B\subset \XX$ viewed from $x=(x_1,x_2,\ldots,x_r)\in\XX\setminus B$ is defined by
\begin{equation}\label{shadowdef}
 \Sh(x:B):=\{z\in\ganz: \ p_i(z)\ne x_i\quad\text{for all }\ i\in\{1,2,\ldots,r\}, \quad  \Ch_{x,z}\cap B\ne\emptyset  \}.
 \end{equation}
It consists of the closure in $\ganz$ of all  Weyl chambers with apex $x$ which intersect $B$ non-trivially. Notice that in view of (\ref{Weylsingular}) we have 
\begin{align}\label{shadowsingular}
 & \Sh(x:B)\cap \rand_\theta = \{\tilde \xi\in\rand_\theta : \sigma_{x_i,\xi_i}(\RR_{\ge 0})\cap p_i(B)\ne\emptyset \quad\text{for all }\ i\in I^+(\theta)\}.
 \end{align}

We next fix a  base point $\xo=(\xo_1,\xo_2,\ldots,\xo_r)\in\XX$. 
For $x\in\XX$ and $t>0$ we denote by $B_x(t)$ the open ball of
radius $t$ centered at $x$. If $h\in\is(\XX_1)\times\is(\XX_2)\times\cdots\times\is(\XX_r)$ is \st all projections $h_i\in\is(\XX_i)$ are axial with translation length $l_i(h)>0$, then $h$ is an axial isometry of the product $\XX=\XX_1\times \XX_2\times\cdots\times \XX_r$ with translation length $l(h)=\sqrt{l_1(h)^2+l_2(h)^2+\cdots+l_r(h)^2}$, and we denote $\widetilde{h^+}, \widetilde{h^-}\in\rand$ its 
attractive and repulsive fixed points. 
If for $i\in\{1,2,\ldots,r\}$ $\ h_i^+, h_i^-\in\rand_i$ denote the attractive and repulsive  
fixed points of the projection $h_i$, then, since for any 
point $x\in \Ax(h)$ and all $n\in \ZZ\setminus\{0\}$ 
\[\thet(x,h^n x)=\Big(\frac{l_1(h)}{l(h)},\frac{l_2(h)}{l(h)},\ldots,\frac{l_r(h)}{l(h)}\Big)=:\widehat L(h) \in E^+,\]  we get 
\[ \widetilde{h^\pm}=(h_1^\pm,h_2^\pm, \ldots,h_r^\pm, \widehat L(h)\big)\in\regrand .\]
So for $h^\pm:=\pi^F(\widetilde{h^\pm})$ we have $h^{\pm}=(h_1^\pm,h_2^\pm,\ldots,h_r^\pm)$. 


The following proposition states that all  Weyl chamber shadows of sufficiently large balls are large in the sense that they contain an open set in $\rand$. This will be crucial in the proof of the shadow lemma. 
Notice that our idea of proof -- which uses Proposition~4.1 of \cite{1308.5584} as a key ingredient -- also considerably simplifies the proof of the analogous statement for one factor (see \cite[Proposition~3.6]{MR1465601} and \cite[Lemma 3.5]{MR2290453}). 
\begin{prp}\label{largeshadows}
Assume that $g=(g_1,g_2,\ldots,g_r)$ and $h=(h_1,h_2,\ldots,h_r) $ are axial isometries of  $\is(\XX_1)\times\is(\XX_2)\times\cdots\times\is(\XX_r)$  
 \st $g_i$ and $h_i$ are independent rank one elements in $\is(\XX_i)$ for all $i\in\{1,2,\ldots,r\}$. Then there exist open 
neighborhoods $U_i\subset\rand_i$ of $h_i^+$, $i\in\{1,2,\ldots,r\}$, a finite set $\Lambda$ in the group $\langle g,h\rangle$ 
generated by $g$, $h$ and $c_0>0$ with the following properties: 

If $U:= U_1\times U_2\times\cdots\times U_r\times E^+\subset \regrand$ and $t\ge c_0$ then for all $y\in\XX\setminus B_\xo(t)$ there exists $\lambda\in\Lambda$ \st
\[  \lambda U\subset \mathrm{Sh}(y:B_\xo(t)).\]
Moreover, if $\theta\in E$ and $U_\theta$ denotes the Cartesian product of the sets $U_i$ with $i\in I^+(\theta)$, then 
\[ \lambda (U_\theta\times\{\theta\}) \subset \mathrm{Sh}(y:B_\xo(t))\cap \rand_\theta.\]
\end{prp}
\prf. \ \ For $i=1,2,\ldots, r$ and $\eta_i\in\{g_i^-,g_i^{+},h_i^-,h_i^{+}\}$ we let $U_i(\eta_i)\subset\ganz_i$ be an arbitrary sufficiently small neighborhood of $\eta_i^+\in\rand_i$ with $\xo_i\notin U_i(\eta_i)$ \st 
all $U_i(\eta_i)$ are pairwise disjoint in $\ganz_i$. Upon taking smaller neighborhoods if necessary Lemma~\ref{joinrankone} provides a constant  $c>0$  \st for every $i\in\{1,2,\ldots,r\}$ any pair of points  in 
distinct sets $U_i(\eta_i)$  can be joined by a rank one geodesic $\sigma_i\subset\XX_i$ with $d(\xo_i,\sigma_i)\le c$. Moreover, according to Lemma~\ref{dynrankone} (b) there exists  a constant   
$N\in\NN$ \st for all  
$i\in\{1,2,\ldots,r\}$
\begin{equation}\label{factorpingpong}  
g_i^{\pm N}\big(\ganz_i\setminus  U_i(g_i^{\mp})\big)\subset U_i(g_i^\pm),\qquad h_i^{\pm N}\big(\ganz_i\setminus  U_i(h_i^{\mp})\big)\subset U_i(h_i^\pm).
\end{equation}
We use induction on $r$ to show the existence of a finite set $\Lambda\subset\langle g,h\rangle$ \st for any $y\in\XX$ one can find $\lambda\in\Lambda$ with 
\[ \lambda y\in U_1(h_1^-)\times U_2(h_2^-)\times\cdots\times U_r(h_r^-).\] 
For $r=1$ we let $y=y_1\in\ganz_1=U_1(h_1^+)\cup \ganz_1\setminus U_1(h_1^+)$ be arbitrary.  If $y_1\in  U_1(h_1^+)$, then from  $U_1(h_1^+)\subset \ganz_1\setminus U_1(g_1^+)$ and~(\ref{factorpingpong}) we get $g_1^{-N} y_1\in U_1(g_1^-)\subset \ganz_1\setminus  U_1(h_1^+)$, hence again by~(\ref{factorpingpong})
\[h_1^{-N} g_1^{-N} y_1\in U_1(h_1^-).\]
If $y_1\in\ganz_1\setminus U_1(h_1^-)$, then~(\ref{factorpingpong}) directly gives $h_1^{-N}y_1\in U_1(h_1^-)$. 
So for $r=1$ the set $\Lambda_1:=\{h^{-N}g^{-N}, h^{-N}\}\subset\langle g,h\rangle$ is the desired finite set.  

Now assume the assertion holds for $r-1$;  we claim that it also holds when $r$ factors are involved.  By the induction hypothesis there exists a finite set \[\Lambda_{r-1}\subset \langle (g_1,g_2,\ldots, g_{r-1}),(h_1,h_2,\ldots, h_{r-1})\rangle<\is(\XX_1)\times\is(\XX_2)\times\cdots\times\is(\XX_{r-1})\] 
\st for all points $(y_1,y_2,\ldots, y_{r-1})\in  \XX_1\times \XX_2\times\cdots\times \XX_{r-1}$
there exists\break  $\lambda'=(\lambda_1',\lambda_2',\ldots,\lambda_{r-1}')\in\Lambda_{r-1}$ \st  $\lambda_i'  y_i \subset U_i(h_i^-)$ for all $i\in\{1,2,\ldots,r-1\}$.
We denote by  $\Lambda'\subset \langle g,h\rangle^+$ the finite set of the same words as in $\Lambda_{r-1}$, but now considered as elements in $\is(\XX_1)\times\is(\XX_2)\times\cdots\times\is(\XX_{r})$, and fix an arbitrary point  $y=(y_1,y_2,\ldots, y_r)\in \XX_1\times \XX_2\times\cdots\times \XX_r$. 
By the properties of $\Lambda_{r-1}$ we know that there exists $\lambda'=(\lambda_1',\lambda_2',\ldots,\lambda_{r-1}',\lambda_r') \in\Lambda'$ \st  
\[ \lambda_i'  y_i \subset U_i(h_i^-)\quad\text{for all} \quad i\in\{1,2,\ldots,  r-1\},\]
but we do not know the position of $\lambda_r' y_r\in\XX_r=U_r(h_r^+)\cup \ganz_r\setminus U_r(h_r^+)$. 

However, as in the case $r=1$ the north-south dynamics~(\ref{factorpingpong}) implies 
\[ h_r^{-N} g_r^{-N} \lambda_r'y_r\in U_r(h_r^-)\quad\text{ or}\quad h_r^{-N} \lambda_r'y_r\in U_r(h_r^-)\]
according to the two cases $\lambda_r'y_r\in U_r(h_r^+)$ or $\lambda_r' y_r\in\ganz_r\setminus U_r(h_r^+)$. Since for $1\le i\le r-1$ 
we have
\[ h_i^{-N} g_i^{-N}\cdot  U_i(h_i^-)\subset U_i(h_i^-) \quad\text{ and}\quad h_i^{-N} \cdot U_i(h_i^-)\subset  U_i(h_i^-)\] 
we conclude that 
 the set $\Lambda_r$ consisting of all words in $g^{-N}, h^{-N}$ of the form
$h^{-N}\lambda'$ or $h^{-N}g^{-N}\lambda'$ with $\lambda'\in\Lambda'$ works. 

So we have shown the existence of  a finite set $\Lambda\subset \langle g,h\rangle$ \st for any $y\in \XX$ there exists  $\lambda=(\lambda_1,\lambda_2,\ldots,\lambda_r)\in\Lambda$ \st  
\[ \lambda^{-1} y \in U_1(h_1^-)\times U_2(h_2^-)\times\cdots\times U_r(h_r^-).\]
In particular, by our choice of the neighborhoods $U_i(h_i^\pm)$, $i=1,2,\ldots r$, every point $z=(z_1,z_2,\ldots,z_r)\in U_1(h_1^+)\times U_2(h_2^+)\times\cdots\times U_r(h_r^+)\subset\ganz_1\times\ganz_2\times\cdots\times\ganz_r$ 
satisfies
\[ d_i(\sigma_{\lambda_i^{-1}y_i, z_i},\xo_i)\le c,\quad\text{ for all}\quad  i\in\{1,2,\ldots,r\}.\]
We next set $d:=\max\{d_i(\xo_i,\lambda_i\xo_i): i\in\{1,2,\ldots,r\}, \lambda=(\lambda_1,\lambda_2,\ldots,\lambda_r)\in\Lambda \}$. Then for $i\in\{1,2,\ldots,r\}$ we have 
\begin{align}\label{inWeylchamber}
d_i(\sigma_{y_i,\lambda_i  z_i},\xo_i)&\le  d_i(\lambda_i \sigma_{\lambda_i^{-1}y_i, z_i},\lambda_i\xo_i)+d_i(\lambda_i\xo_i,\xo_i)\nonumber\\
&<  d_i(\sigma_{\lambda_i^{-1}y_i, z_i},\xo_i) + d \le c+d.
\end{align}
We set  $c_0:=\sqrt{r}(c+d)$ and $U_i:=U_i(h_i^+)\cap\rand_i$ for  $i\in\{1,2,\ldots,r\}$. If $\theta\in E$ and $\tilde\zeta\in U_\theta\times\{\theta\}\subset\rand_\theta$, then according to (\ref{inWeylchamber}) its projections $\zeta_i\in U_i$, $i\in I^+(\theta)$, satisfy 
\[ d_i(\sigma_{y_i,\lambda_i  \zeta_i},\xo_i)\le c+d,\]
hence by definition~(\ref{Weylsingular}) of the Weyl chamber with apex $y$ in the class of $\lambda\tilde\zeta$ we conclude that for all $t\ge c_0$ 
\[ \Ch_{y,\lambda \tilde\zeta}\cap B_\xo(t)\ne\emptyset,\quad\text{ and hence}\quad 
\lambda\tilde\zeta\in\Sh(y:B_{\xo}(t))\cap\rand_\theta.\]
The claim for $U=U_1\times U_2\times\cdots\times U_r\times E^+ \subset\regrand$ follows from the fact that 
\[ U=\bigcup_{\theta\in E^+} U_\theta.\tag*{$\square$}\] 

For $\theta=(\theta_1,\theta_2,\ldots,\theta_r)\in E=\{\theta\in\RR^r_{\ge 0}:\Vert \theta\Vert=1\}$ we recall from~(\ref{randtheta}) the definition of the set $\rand_\theta\subset\rand$ which is homeomorphic to the Cartesian product of the geometric boundaries $\rand_i$ with  $i\in I^+(\theta)=\{i\in\{1,2,\ldots,r\}: \theta_i>0\}$. 
The following easy lemma relates the Busemann function~(\ref{buseman}) of the product to the Busemann functions on the factors. 
We include a proof for the convenience of the reader.
\begin{lem}\label{RelProdFactorr}  
Let $\theta=(\theta_1,\theta_2,\ldots,\theta_r)\in E$,  $x=(x_1,x_2,\ldots,x_r)$, $y=(y_1,y_2,\ldots,y_r)\in\XX$ and $\tilde\xi\in\rand_\theta$. 
If $\xi_i$ denotes the projection of $\tilde\xi$ to $\rand_i$ then 
\[\bs_{\tilde\xi}(x,y)=\sum_{i\in I^+(\theta)} \theta_i \cdot \bs_{\xi_i}(x_i,y_i).\]
\end{lem}
\prf. \ \ Notice that from the definition of the Busemann functions in $\XX_i$, $i\in I^+(\theta)$, 
we have
\begin{align*}
\bs_{\xi_i}(x_i,y_i)&=\lim_{s\to\infty} \big( s \theta_i-d_i(y_i,\sigma_{x_i,\xi_i}(s\theta_i))\big).
\end{align*}
For convenience we set $I^0(\theta)=\{1,2,\ldots, r\}\setminus I^+(\theta)$. Since  $\displaystyle \sum_{i\in I^+(\theta)}\theta_i^2= \sum_{i=1}^r\theta_i^2= 1 $  we get 
\vspace*{-4mm}
 \begin{align*}  
\hspace*{6mm} \big(s-d(y,\sigma_{x,\tilde\xi}(s))\big) & \big(s+d(y,\sigma_{x,\tilde\xi}(s))\big) =s^2-d(y,\sigma_{x,\tilde\xi}(s))^2\\
& =  \displaystyle s^2\sum_{i\in I^+(\theta)} \theta_i^2 -\sum_{i\in I^+(\theta)} d_i(y_i,\sigma_{x_i,\xi_i}(s\theta_i))^2-\sum_{i\in I^0(\theta)} d_i(y_i,x_i)\\
&=\sum_{i\in I^+(\theta)} s^2\theta_i^2-d_i(y_i,\sigma_{x_i,\xi_i}(s\theta_i))^2-\sum_{i\in I^0(\theta)} d_i(y_i,x_i) .
\end{align*}
So the assertion is proved if we show that for all $i\in I^+(\theta)$
 \begin{align*} &  \lim_{s\to\infty}\frac{s\theta_i+d_i(y_i,\sigma_{x_i,\xi_i}(s\theta_i))}{s+d(y,\sigma_{x,\tilde\xi}(s))}=\theta_i;
 \end{align*}
 this claim follows immediately from the triangle inequalities
\begin{align*} 
&s\theta_i-d_i(y_i,x_i)\le d_i(y_i,\sigma_{x_i,\xi_i}(s\theta_i))\le s\theta_i +d_i(y_i,x_i),\\
& s-d(y,x)\le d(y,\sigma_{x,\tilde\xi}(s)) \le s+d(y,x).\tag*{$\square$}\\[-4mm]
\end{align*}

To simplify notation in the sequel we further define for 
$x=(x_1,x_2,\ldots,x_r)$,\break $y=(y_1,y_2,\ldots,y_r)\in\XX$ and $\tilde\xi\in\rand_\theta$  the 
{\hl Busemann vector} 
\begin{equation}\label{busvector}\bv_{\tilde\xi} (x,y)\end{equation} 
as the unique vector in $\RR^r$ with 
$i$-th coordinate equal to $\bs_{\xi_i}(x_i,y_i)$ for $i\in I^+(\theta)$, and $i$-th coordinate equal to zero for all $i\in I^0(\theta)$. 

Notice that for $\tilde\xi\in\regrand$, the Busemann vector $\bv_{\tilde\xi}$ is independent of the slope of $\tilde\xi$; it only depends on $(\xi_1,\xi_2,\ldots,\xi_r)=\piF(\tilde\xi)\in\Frand$. Moreover, by the cocycle identity 
for the Busemann function we get
\[ \bv_{\tilde\xi}(x, z)=\bv_{\tilde\xi}(x, y)+\bv_{\tilde\xi}(y,z)\quad\text{ for all}\quad x,y,z\in\XX.\]
We also remark that if $\langle\cdot,\cdot\rangle$ denotes the Euclidean inner product of $\RR^r$, then the formula in Lemma~\ref{RelProdFactorr} can be rewritten as
\begin{equation}\label{RelProdFactor}
 \bs_{\tilde\xi}(x,y) =\langle \bv_{\tilde\xi}(x,y), \theta\rangle.
\end{equation}
In the sequel we will also need the following 
\begin{df}\label{dirdist} The {\hl directional distance} of the ordered pair $(x, y)\in
\XX\times \XX$ \wrt the slope $\theta\in E$ is defined by
\[\begin{array}{rclrl}
 \bs_{\theta}:\quad\XX\times \XX &\to &\RR && \\[1mm]
 (x, y) &\mapsto & \bs_{\theta} (x, y)\;& := &\, \langle H(x,y), \displaystyle \theta\rangle =\sum_{i=1}^r \theta_i\cdot   d_i(p_i(x),p_i(y)).
\end{array}\] 
\end{df}
In particular,  if $\theta\in E$ has $i$-th coordinate $1$ and all other coordinates zero then $\bs_\theta(x,y)=d_i(p_i(x),p_i(y))$.

By $\big(\is(\XX_1)\times\is(\XX_2)\times\cdots\times\is(\XX_r)\big)$-invariance of the distance vector we immediately get that 
$$ \bs_\theta(gx,gy)=\bs_\theta(x,y)$$ for all $x,y\in\XX$ and $g\in  \is(\XX_1)\times\is(\XX_2)\times\cdots\times\is(\XX_r)$. Moreover, the symmetry and triangle inequality for the distances $d_1,d_2,\ldots, d_r$ 
directly imply the symmetry and triangle inequality for $\bs_\theta$.
The following important proposition states that for $\theta\in E^+=\{\theta\in\RR^r_{>0}: \Vert \theta\Vert=1\}$ the directional distance $\bs_\theta$ is in fact a distance.
\begin{prp}\label{dirisdistance}
For  $\theta\in E^+$ the directional distance $\bs_\theta$ is a distance.
\end{prp}
\prf. \  \ Let $x=(x_1,x_2,\ldots,x_r)$, $y=(y_1,y_2,\ldots,y_r)\in\XX$. We clearly have 
\[ \bs_\theta(x,y)=\sum_{i=1}^r \theta_i \cdot d_i(x_i,y_i)\ge 0,\] 
because all terms involved are non-negative. 
Moreover, if $\bs_\theta(x,y)=0$, then $\theta_i>0$ for all $i\in\{1,2,\ldots,r\}$ imply $ d_1(x_1,y_1)=d_2(x_2,y_2)=\cdots=d_r(x_r,y_r)=0$, hence  $x=y$.

Finally, we have already noticed that the symmetry and triangle inequality follow directly from the symmetry and triangle inequality for the distances $d_i$. 
\qed\\[-1mm]

The following easy facts will be convenient in the sequel.
\begin{lem} \label{dirBus}
Let $x,y \in\XX$  and $\tilde\xi\in \rand_\theta$ for some $\theta\in E$. Then 
$$ y\in \Ch_{x,\tilde\xi}
\ \iff\quad \bs_\theta(x,y)=\bs_{\tilde\xi}(x,y).$$
\end{lem}
\prf.\ \  We write $x=(x_1,x_2,\ldots,x_r)$, $y=(y_1,y_2,\ldots,y_r)$ and $\theta=(\theta_1,\theta_2,\ldots,\theta_r)$. Recall that $I^+(\theta)=\{i\in\{1,2,\ldots,r\}: \theta_i>0\}$, and that for $i\in I^+(\theta) $ $\xi_i$ denotes the projection of $\tilde\xi$ to $\rand_i$.  
Lemma~\ref{RelProdFactorr}  and the estimates $\bs_{\xi_i}(x_i,y_i)\le d_i(x_i,y_i)$ 
imply 
\begin{equation}\label{ineqbusedir}
\bs_{\tilde\xi}(x,y)=\sum_{i\in I^+(\theta)}\theta_i\cdot \bs_{\xi_i}(x_i,y_i) \le \sum_{i=1}^r \theta_i\cdot d_i(x_i,y_i) =\langle  H(x,y),\theta\rangle=\bs_\theta(x,y).
\end{equation} 
So we have equality in~(\ref{ineqbusedir}) if and only if for all $i\in I^+(\theta)$ the equality \[\bs_{\xi_i}(x_i,y_i)=d_i(x_i,y_i)\] holds; this is precisely the case when $y_i$ is a point on the geodesic ray $\sigma_{x_i\xi_i}$ in $\XX_i$. 
Therefore -- by definition~(\ref{Weylsingular}) of the Weyl chamber with apex $x$ in the class of $\tilde\xi$ --  equality in~(\ref{ineqbusedir}) is equivalent to $y\in  \Ch_{x,\tilde\xi}$.  \qed\\[-1mm]

If some of the factors $\XX_i$ are geodesically complete, the previous lemma allows to give the following nice geometric interpretation of the directional distance. 
\begin{cor}\label{dirdistbusemann}
Fix $\theta\in E$ and assume that  $\XX_i$ is geodesically complete for all\break $i\in I^+(\theta)$. Then for all $x,y\in \XX$  we have
$$\bs_\theta(x,y)= \max \{ \bs_{\tilde \xi} (x,y): \tilde\xi\in\rand_\theta\}.$$
\end{cor}
\prf.\ \ We first fix $i\in I^+(\theta)$. Since $\XX_i$ is geodesically complete, every point\break $y_i\in\XX_i\setminus \{x_i\}$ belongs to a geodesic ray $\sigma_{x_i,\xi_i}$ with $\xi_i\in \rand_i$ the unique extension $\sigma_{x_i,y_i}(\infty)$ of the geodesic in $\XX_i$ joining $x_i$ to $y_i$.  If $y_i=x_i$ one  may choose an arbitrary point $\xi_i\in\rand_i$. 

In this way  every $y\in\XX$ determines a (not necessarily unique) boundary point $\tilde\xi\in\rand_\theta$ with projections $\xi_i\in\rand_i$, $i\in I^+(\theta)$; by choice of $\tilde\xi$ and definition~(\ref{Weylsingular}) we clearly have $y\in \Ch_{x,\tilde\xi}$ and hence, by Lemma~\ref{dirBus},
\[\bs_\theta(x,y)=\bs_{\tilde\xi}(x,y).\] 
Inequality~(\ref{ineqbusedir}) then proves the claim. \qed\\[-1mm]

Recall the definition of Weyl chamber shadows from~(\ref{shadowdef}) and~(\ref{shadowsingular}). The following lemma will be needed in the proof of the shadow lemma Theorem~\ref{shadowlemma}.
\begin{lem}\label{esti}
Let $c>0$, $z=(z_1,z_2,\ldots,z_r) \in\XX$ \st  $d(o,z)>c$, $\theta\in E$ and\break $\tilde\eta\in\Sh\big(o:B_{z}(c)\big)\cap\rand_\theta$ with projections $\eta_i\in\rand_i$, $i\in I^+(\theta)$. Then we have
\[   0\le d_i(\xo_i,z_i)-\bs_{\eta_i}(\xo_i,z_i)<2c\quad\text{ for all} \quad i\in I^+(\theta).\]
\end{lem} 
\prf.\ \  By definition $\tilde\eta\in \Sh\big(o:B_{z}(c)\big)$ if and only if $ \Ch_{\xo,\tilde\eta}\cap B_z(c)\ne\emptyset$. Hence if $\tilde\eta\in\rand_\theta$  then for all $i\in I^+(\theta)$ there exists
$t_i\ge 0$ \st $y_i:=\sigma_{\xo_i,\eta_i}(t_i)\in p_i(B_z(c))$. Necessarily we have
\[ d_i(z_i,\sigma_{\xo_i,\eta_i})\le d_i(z_i,y_i)<c\quad\text{ for all} \quad i\in I^+(\theta),\]
hence the claim follows from Lemma~\ref{distbusbounded}.\qed

\section{The exponent of growth}\label{ExpGrowth}

For the remainder of the article  $\XX$ will be a product of locally compact 
Hadamard spaces $\XX_1$, $\XX_2,\ldots, \XX_r$,  and $\Gamma<\is(\XX_1)\times\is(\XX_2)\times\cdots\times \is(\XX_r)$ a group acting properly discontinuously by isometries on $\XX$ which contains two elements $h=(h_1,h_2,\ldots,h_r)$ and $g=(g_1,g_2,\ldots, g_r)$ \st for $i\in\{1,2,\ldots,r\}$ $g_i$ and $h_i$ are independent rank one isometries of $\XX_i$. We further fix a base point $\xo=(\xo_1, \xo_2,\ldots,\xo_r)\in\Ax(h) $. 

We recall that the geometric limit set of  $\Gamma$ 
is defined by
$\Lim:=\overline{\Gamma\at x}\cap\rand$, where $x\in \XX$ is arbitrary.
In this section we recall the notion of  exponent of  growth introduced in \cite{MR2629900} and \cite{1308.5584} and give an important criterion for divergence or convergence of certain sums over $\Gamma$. 
This will play a central role in the construction of (generalized) 
conformal densities 
in Sections~\ref{ClasPatSulconst} and \ref{genPseries}.

We recall the notation introduced in Section~\ref{prodHadspaces}; in particular, we denote $E\subset \RR^r$ the set of unit vectors in $\RR_{\ge 0}^r$.  For $x,y\in\XX$, $\theta \in E$ and  $\eps>0$ we first set
\[ \Gamma(x,y;\theta,\eps):=\{\gamma\in\Gamma: \gamma y\ne x\quad \mbox{and}\ \ \Vert \thet(x,\gamma y)-\theta \Vert <\eps \}.\]

In order to define the exponent of growth of $\Gamma$ of slope $\theta$ we set
$$ \delta_\theta^\eps(x,y):=\inf \{s>0: \sum_{\gamma\in
  \Gamma(x,y;\theta,\eps)}
  \e^{-s d(x,\gamma y)} \ \ \mbox{converges}\}.$$
If $\delta(\Gamma)$ denotes the {\hd critical exponent} of $\Gamma$ defined by
\begin{equation}\label{critexp}
  \delta(\Gamma):=\inf\{s>0: \sum_{\gamma\in\Gamma}\e^{-sd(\xo,\gamma\xo)}\ \mbox{converges}\},
\end{equation}
we clearly have $\delta_\theta^\eps(x,y)\le \delta(\Gamma)$ with equality if $\eps>\sqrt2$.
Moreover, an easy calculation 
shows that  $\delta_\theta^\eps(x,y)$ is related to the numbers  
\[\Delta N_\theta^\eps(x,y;n):=\#\{ \gamma\in\Gamma:\ n-1\le d(x,\gamma y)<n,\ \Vert\thet(x,\gamma y)-\theta\Vert<\eps\} \]
with $n\in\NN$, $n\ge 2$ via 
\begin{equation}\label{defbylimsup}
\delta_{\theta}^{\eps}(x,y)=\limsup_{n\to\infty}\frac{\ln
  \Delta N_{\theta}^{\eps}(x,y;n)}{n}.
  \end{equation}
Recall that the {\hl exponent of growth of $\Gamma$ of slope $\theta$} is  defined by 
\[\delta_{\theta}(\Gamma):=\lim_{\eps\to
  0}\delta_{\theta}^{\eps}(\xo,\xo).\]
 Notice that this number  $\delta_\theta(\Gamma)$  does not depend on the choice of arguments 
of $\delta_{\theta}^\eps$ by elementary geometric estimates; 
it can be interpreted as 
an  exponential growth rate of the number of orbit points which are {``}close" to a geodesic ray in the class of a boundary point with slope $\theta$. 
Moreover, we clearly have 
$\delta_\theta(\Gamma)\le \delta(\Gamma)$ for all $\theta\in E$. 

Furthermore, we recall the following properties from Section~7 in \cite{1308.5584}:\\[3mm]
\begin{pros}
\begin{enumerate}
\item[(a)] $\Lim\cap  \rand_\theta\ne \emptyset$ if and only if $\delta_{\theta}(\Gamma)\ge 0$.
\item[(b)] The map $E\to\RR,$ $\hat\theta\mapsto \delta_{\hat\theta}(\Gamma)$
is upper semi-continuous.
\end{enumerate}
\end{pros}

It will turn out useful to consider the homogeneous extension  $\Psi_\Gamma:\RR_{\ge 0}^r\to \RR$ of the map $E\to\RR,$ $\hat\theta\mapsto \delta_{\hat\theta}(\Gamma)$.   
Theorem~7.6 in \cite{1308.5584} states that $\Psi_\Gamma$ is concave. This implies in particular that there exists a unique $\theta^*\in E$ such that 
$\delta_{\theta^*}(\Gamma)=\max\{\delta_\theta(\Gamma):\theta\in E\}$. The following important proposition will play a key role in the proof of Theorem~A and for the construction of generalized conformal densities.
Recall definitions~(\ref{distancevector}) and~(\ref{slopeangle}) for the distance vector and the direction of a pair of points in $\XX$. 
\begin{prp}\label{convdiv}
Let  $f: \RR_{\ge 0}^r\to \RR $ be a continuous homogeneous function, $D\subset E$ a  a relatively  open set and put
$\,\Gamma_D:=\{\gamma\in\Gamma :\  \gamma\xo\ne \xo,\  \thet(\xo,\gamma\xo)\in D\}$.
\begin{enumerate}
\item[(a)]
If there exists $\hat\theta\in D$ \st $f(\hat\theta)<\delta_{\hat\theta}(\Gamma),$
then the sum $\displaystyle\sum_{\gamma\in\Gamma_D}\e^{-f(H(\xo,\gamma\xo))}$ diverges.
\item[(b)]
If $f(\theta)> \delta_{\theta}(\Gamma)$ 
for all $\theta\in \overline{D}$,
 then the sum  $\displaystyle\sum_{\gamma\in\Gamma_D}\e^{-f(H(\xo,\gamma\xo))}$ converges.
\end{enumerate} 
\end{prp}
\prf.\ \  For $\gamma\in\Gamma$ we abbreviate $\displaystyle \theta_\gamma:=\thet(\xo,\gamma\xo)=\frac{H(\xo,\gamma\xo)}{d(\xo,\gamma\xo)}$.
\begin{enumerate}
\item[(a)] Let $\hat\theta \in D$  \st 
$\  f(\hat\theta)<\delta_{\hat\theta}(\Gamma).$  
Since $\displaystyle\delta_{\hat\theta}(\Gamma)=\lim_{\eps\to
  0}\delta_{\hat\theta}^\eps(o,o)$, there exists
$\eps>0$ and $ \hat s\in\RR\ $ \st for 
$\gamma\in\Gamma_D\ $ with  
$\ \Vert \theta_\gamma-\hat\theta \Vert <\eps\ $  we have  
$$  f(\theta_\gamma)< \hat s<\delta_{\hat\theta}^\eps(\xo,\xo).$$
Since $f(H(\xo,\gamma\xo))= f(\theta_\gamma) \cdot d(\xo,\gamma\xo)$ we estimate   
\[ \sum_{\gamma\in\Gamma_D}\e^{-f(H(\xo,\gamma\xo))}
> 
\sum_{\gamma\in\Gamma(\xo,\xo;\hat\theta,\eps)}\e^{- \hat s d(o,\gamma o)},\]
and the latter sum diverges because $ \hat s<\delta_{\hat\theta}^\eps(o,o)$.

\item[(b)]
Let $ \hat\theta\in\overline{D}$. Since $\displaystyle f(\hat\theta)> \delta_{\hat\theta}(\Gamma)=\lim_{ \eps\to 0}\delta_{\hat\theta}^ \eps(\xo,\xo)$, there exists  $ \eps'>0$ and
  $ \hat s<f(\hat\theta)$ \st 
\begin{equation}\label{sandwich} 
\delta_{  \hat\theta}^{ \eps'}(o,o)< \hat s< f(\hat\theta).
\end{equation}
For $\theta\in E$ and $\eps>0$ we set $B_{\theta}( \eps):=\{\theta'\in E : \Vert\theta'-\theta\Vert< \eps\}$. The continuity of the function $f$
and inequality~(\ref{sandwich}) imply
the existence of $ \hat\eps< \eps'$ \st for any $\theta\in
B_{  \hat\theta}( \hat\eps)$  we have $ \hat s<f(\theta)$. Hence for all $z\in\XX$ with $\theta_z:=\thet(\xo,z)\in B_{  \hat\theta}( \hat\eps)$ we have
$$f(\theta_z)= \frac{f(H(\xo,z))}{d(\xo,z)}>\hat s> \delta_{  \hat\theta}^{ \eps'}(o,o) \ge \delta_{  \hat\theta}^{ \hat \eps}(o,o).$$
We now choose a sequence $(\theta_j)\subset\overline{D}$ and corresponding
sequences  $( \eps_j)\subset\RR_{>0}$ and $(s_j)\subset \RR_{>0}$ \st for every 
$\theta\in B_{ \theta_j}( \eps_j)$   we have 
\[\delta_{ \theta_j}^{ \eps_j}(o,o)<s_j<f(\theta)\]
and such that the sets $B_{\theta_j}(\eps_j)$, $j\in\NN$, cover $\displaystyle \overline{D}$. 
By compactness of $\overline D$ there exists a finite set $J\subset\NN$ with 
$\displaystyle  \overline{D}\subset \bigcup_{j\in J} B_{ \theta_j}( \eps_j)$, 
and we conclude 
\begin{align*}
\sum_{\gamma\in\Gamma_D}\e^{-f(H(\xo,\gamma\xo))}& \le
\sum_{j\in J}
\sum_{\gamma\in\Gamma(\xo,\xo;\theta_j,\eps_j)}
\e^{-f(H(\xo,\gamma\xo))}\\
&\le\sum_{j\in J} \sum_{\gamma\in\Gamma(\xo,\xo;\theta_j,\eps_j)}
\e^{-s_j d(o,\gamma o)} <\infty,
\end{align*}
because $s_j>\delta_{ \theta_j}^{ \eps_j}(o,o)$ for $j\in J$.\qed\\[-1mm]
\end{enumerate}
Taking $D=E$ and $f(H)=s\cdot \Vert H\Vert\ $ we obtain as a corollary that 
\[ \delta(\Gamma)=\max\{\delta_\theta(\Gamma):\theta\in E\}=\delta_{\theta^*}(\Gamma).\] 
 We conclude this section with two illustrative examples. \\[-1mm]

\noindent{\sc Example 1} (see \cite[Section 7]{1308.5584}) $\quad$
We let $\XX=\XX_1\times \XX_2\times\cdots\times\XX_r$ be a product of Hadamard manifolds with pinched negative curvature, and assume that for all\break $i\in\{1,2,\ldots,r\}$ a discrete convex cocompact group  $\Gamma_i<\is(\XX_i)$   with critical exponent $\delta_i>0$ is given. Then the exponent of growth of slope $\theta\in E$ for  the product group $\Gamma:=\Gamma_1\times\Gamma_2\times\cdots\times\Gamma_r$  satisfies
\[  \delta_\theta(\Gamma)=\sum_{i=1}^r \delta_i\theta_i .\]
Using Lagrange multipliers one can easily show that this number is maximal for\break $\theta^*\in E^+$ with coordinates
\[ \theta_i^*=\frac{\delta_i}{\delta_1^2+\delta_2^2+\cdots+\delta_r^2},\quad i\in\{1,2,\ldots,r\};\]
in particular we have 
\[\delta(\Gamma)=\delta_{\theta^*}(\Gamma)=\sqrt{\delta_1^2+\delta_2^2+\cdots+\delta_r^2}.\] 
The homogeneous function $\Psi_\Gamma:\RR_{\ge 0}^r\to \RR$ is simply the linear functional defined by taking the inner product $\langle\cdot,\cdot\rangle$ in $\RR^r$ with the unique vector in $\RR^r_{>0}$ with coordinates $\delta_i$.  \\[2mm]

\noindent{\sc Example 2}  $\quad$  Consider a product of hyperbolic planes $\XX=\HH^2\times \HH^2$ and a Hilbert modular group $\Gamma<\is(\XX)$. Then $\Gamma$ is an irreducible 
non-uniform lattice in a higher rank symmetric space,
hence from Proposition 7.2 and 7.3 in \cite{MR1675889} we know that 
\[\Psi_\Gamma=\langle \left(\ba{c}1\\1\ea\right),\cdot\rangle,\quad\text{so}\quad  
\delta_\theta(\Gamma)=\theta_1 +\theta_2.\] 
Here $\delta_\theta(\Gamma)$ is maximal for 
\[ \theta^*=\frac{1}{\sqrt{2}}  \left(\ba{c} 1\\1\ea\right) \in E^+, \] so we get $\delta(\Gamma)=\delta_{\theta^*}(\Gamma)=\sqrt2$.

\section{The classical Patterson-Sullivan construction}\label{ClasPatSulconst}

In this section we will  construct  a conformal density for $\Gamma$ using an idea originally due to  S.~J.~Patterson (\cite{MR0450547}) in the context 
of Fuchsian groups. Taking advantage of Proposition~\ref{convdiv} we will be able to describe precisely its support and hence 
 prove Theorem~A. 

Recall that a $\Gamma$-invariant conformal density of dimension $\delta\ge 0$ is  a
map $\mu$ from $\XX$ to the cone $\MM^+(\rand)$ of positive   finite Borel measures
on $\rand$ \st 
$\supp(\mu_{\xo})\subset\Lim$, $\gamma_*\mu_x=\mu_{\gamma x}$ for all $\gamma\in\Gamma$, $x\in\XX$ and
$$ \frac{d\mu_x}{d\mu_\xo}(\tilde\eta)=\e^{\delta \bs_{\tilde \eta}(\xo,x)}\qquad\mbox{for all}\ \tilde\eta\in\supp(\mu_{\xo}),\  x\in\XX.$$

In order to construct a $\Gamma$-invariant conformal density of dimension $\delta(\Gamma)$ we first suppose that we are given a map $b:\Gamma\to\RR$, $\gamma\mapsto b_\gamma$ \st the sum 
\begin{equation}\label{PatSum}\sum_{\gamma\in\Gamma} \e^{-sb_\gamma}
\end{equation}
has exponent of convergence $s=1$ (which means that it converges for $s>1$ and diverges for $s<1$). The following useful lemma states that if the above sum converges for $s=1$, then we can slightly modify it to 
obtain a sum which diverges for $s\le 1$ and converges for $s>1$.   

\begin{lem}[\ (Patterson \cite{MR0450547}\label{Pt})]
If the sum~(\ref{PatSum}) 
has exponent of convergence $s=1$, then there exists a non-decreasing continuous function $h:[0,\infty)\to [1,\infty)$ \st 
\begin{enumerate}
\item[\rm(i)]
$\displaystyle\sum_{\gamma\in \Gamma} \e^{-s b_\gamma}
h(\e^{b_\gamma})$ has exponent of convergence $s=1$ and
diverges at $s=1$;
\item[\rm (ii)]
for any $\alpha>0$ there exists $r_0>0$  \st  for $r\ge r_0$ and $t>1$
 $$h(rt)\le t^\alpha h(r).$$
\end{enumerate}
\end{lem}
Notice that if the sum~(\ref{PatSum}) already diverges at $s=1$, then $h$ can be chosen as the constant function identical to $1$. 

Recall the definition of the exponent of growth of $\Gamma$ and its properties from Section~\ref{ExpGrowth}. We have already noticed that there exists a unique 
unit vector $\theta^*\in E$ \st $\delta(\Gamma)=\delta_{\theta^*}(\Gamma)$. 

Following the original idea of Patterson \cite{MR0450547}, we apply the above lemma to the map 
\begin{equation}\label{bClas}
b:\Gamma\to\RR, \quad \gamma\mapsto \delta(\Gamma)\cdot d(\xo,\gamma\xo). 
\end{equation}
Then by definition~(\ref{critexp}) of the critical exponent $\delta(\Gamma)$ the sum
\[\displaystyle\sum_{\gamma\in\Gamma}\e^{-sb_\gamma}=\sum_{\gamma\in\Gamma}\e^{-s\delta(\Gamma)\cdot d(\xo,\gamma\xo)}\] 
has exponent of convergence $s=1$.  Let  $h:[0,\infty)\to [1,\infty)$ be a non-decreasing function as in Patterson's Lemma  above and define
\begin{equation}\label{Summs} 
\Summ^s:=\sum_{\gamma\in\Gamma}\e^{-sb_\gamma}h(\e^{b_\gamma}).
\end{equation}
If  $D$ denotes the unit Dirac point measure, then for  $s>1$ we get a probability measure on $\XX$ by setting
\begin{equation}\label{probmeasdef}
 \mu_\xo^s:=\frac1{\Summ^s}\sum_{\gamma\in\Gamma}\e^{-sb_\gamma }h(\e^{b_\gamma})D(\gamma\xo).
 \end{equation}
Notice that by construction any  weak accumulation point $\mu_\xo$ of $\mu_\xo^{s}$ as $s\searrow 1$ is a probability measure  on $\rand$ with 
$\supp(\mu_\xo)\subset\Lim$.

Before we continue with the construction of a $\Gamma$-invariant conformal density 
we state an auxiliary lemma which will be useful in the sequel. For a topological space $Y$ we denote $(\Cnt^{0}(Y),\Vert\cdot\Vert_\infty)$ the
space of real
valued continuous functions on $Y$ with norm
$\Vert f\Vert_\infty=\sup\{|f(y)| :  y\in Y\}$,
$f\in \Cnt^{0}(Y)$. 
\begin{lem}\label{supnorm}
Fix $x,y\in \XX$ and $s>0$. Let $h:[0,\infty)\to [1,\infty)$ be a non-decreasing function as in Patterson's Lemma~\ref{Pt}, and $b:\XX\times \XX\to \RR$ a continuous map with the property
\begin{equation}\label{bDiffgen}
 |b(x,z)-b(y,z)|\le C\cdot d(x,y)\quad\text{for all}\quad z\in \XX.
 \end{equation}
Then the continuous function 
\begin{equation}\label{gxysDef}
 g_{x,y}^s:\XX\to\RR,\quad z \mapsto \frac{\e^{-sb(x,z)}h(\e^{b(x,z)})}{\e^{-sb(y,z)}h(\e^{b(y,z)})}
 \end{equation}
satisfies
\[ \lim_{s\to 1} \Vert g_{x,y}^s-g_{x,y}^1\Vert_\infty =0.\]
\end{lem}
\prf.\ \  We first remark that for all $z\in \XX$ we have
\begin{align*} 
\frac{g_{x,y}^s(z)}{g_{x,y}^1(z)} &=\frac{\e^{s(b(y,z)-b(x,z))}}{\e^{b(y,z)-b(x,z)}}=\e^{(s-1)(b(y,z)-b(x,z))},
\end{align*}
so by the property~(\ref{bDiffgen}) of $b$ we get
\begin{align*} 
\e^{-|s-1|C\cdot d(x,y)} &\le \frac{g_{x,y}^s(z)}{g_{x,y}^1(z)} 
\le  \e^{|s-1|C\cdot d(x,y)}.
\end{align*}
Using the inequality $\e^t+\e^{-t}\ge 2$ we obtain for all $z\in \XX$ 
\[ \Big| \frac{g_{x,y}^s(z)}{g_{x,y}^1(z)}-1\Big| \le \e^{|s-1|C\cdot d(x,y)}-1.\]
By Patterson's Lemma (ii) there exists $r_0>0$ \st for all $r\ge r_0$ and $t>1$ we have
\[ h(rt)\le t^2 h(r)\le \e^{2t}h(r).\]
So for $z\in \XX$ \st $b(y,z)\ge \ln(r_0)$ we get
\[ \frac{h(\e^{b(x,z)})}{h(\e^{b(y,z)})}\le \e^{2(b(x,z)-b(y,z))}\le \e^{2C\cdot d(x,y)}.\]
If $z\in\XX$ satisfies $b(y,z)< \ln(r_0)$, then 
\[ b(x,z)=(b(x,z)-b(y,z))+b(y,z)\le C\cdot d(x,y)+\ln(r_0),\]
hence since $h$ is a non-decreasing function $\ge 1$ 
\[ \frac{h(\e^{b(x,z)})}{h(\e^{b(y,z)})}\le h(\e^{b(x,z)})\le h(\e^{C\cdot d(x,y)}r_0) .\]
This implies that there exists a constant $K>1$ (which only depends on $d(x,y)$) \st
for all $z\in\XX$
\[ \frac{h(\e^{b(x,z)})}{h(\e^{b(y,z)})}\le K.\]
We conclude
\begin{align*}
\big| g_{x,y}^s(z) - g_{x,y}^1(z)\big| &= g_{x,y}^1(z)\cdot  \Big| \frac{g_{x,y}^s(z)}{g_{x,y}^1(z)}-1\Big|\le\frac{\e^{-b(x,z)}h(\e^{b(x,z)})}{\e^{-b(y,z)}h(\e^{b(y,z)})}\big( \e^{|s-1|C\cdot d(x,y)}-1\big)\\
&\le \e^{C\cdot d(x,y)}\cdot K \big( \e^{|s-1|C\cdot d(x,y)}-1\big),
\end{align*}
so\\[1mm]
$\hspace*{2cm} \Vert g_{x,y}^s-g_{x,y}^1\Vert_\infty=\displaystyle\sup_{z\in\XX} \big| g_{x,y}^s(z) - g_{x,y}^1(z)\big|\ \longrightarrow\ 0\quad\text{as}\quad s\to 1.$\qed\\[1mm]
In order to obtain a $\Gamma$-invariant conformal density, we imitate the construction~(\ref{probmeasdef}) and define for 
$x,z \in X$ and $\gamma\in\Gamma$
\[ b(x,z):= \delta(\Gamma)\cdot d(x,z),\qquad b_\gamma:=b(\xo,\gamma\xo).\]
Notice that by the triangle inequality for the distance function the map $b:\XX\times\XX\to \RR$ is continuous and satisfies property~(\ref{bDiffgen}) with $C=\delta(\Gamma)$.

For $s>1$, $x\in\XX$ and with $\Summ^s$ as defined in~(\ref{Summs})
we get a family of positive finite Borel measures 
on $\overline{\XX}$ via
\[ \mu_x^s:=\frac1{\Summ^s}\sum_{\gamma\in\Gamma} \e^{-sb(x,\gamma\xo)}h(\e^{b(x,\gamma\xo)})D(\gamma\xo);\]
in particular, $\mu_\xo^s$ is precisely the probability measure defined in~(\ref{probmeasdef}). 
For fixed $s>1$ the measures $\mu_x^s$, $x\in\XX$, are $\Gamma$-equivariant by construction and absolutely
continuous \wrt each other with Radon Nikodym derivative
\begin{equation}\label{RadNikclassical} 
\frac{\d \mu_x^s}{\d \mu_y^s}:\ \supp(\mu_y^s)\to\RR,\qquad z\mapsto \frac{\e^{-sb(x,z)}h(\e^{b(x,z)})}{\e^{-sb(y,z)}h(\e^{b(y,z)})}=g_{x,y}^s(z), 
\end{equation}
where $g_{x,y}^s$ is the continuous function defined by~(\ref{gxysDef}) in Lemma~\ref{supnorm}. 

Moreover, we have the following
\begin{lem}\label{extendscont} For fixed $x,y\in \XX$ and $s>0$  the continuous function 
\[ g_{x,y}^s:\XX\to\RR,\quad z\mapsto \frac{\e^{-sb(x,z)}h(\e^{b(x,z)})}{\e^{-sb(y,z)}h(\e^{b(y,z)})}\]
extends continuously to  $\ganz$.  Moreover, if $(z_n)\subset\XX$ is a sequence converging to $\tilde\eta\in\rand$, then
\[ \lim_{n\to\infty} g_{x,y}^s(z_n)=\e^{-s\delta(\Gamma)\bs_{\tilde\eta}(x,y)}=\e^{s\delta(\Gamma)\bs_{\tilde\eta}(y,x)}.\]
\end{lem}
\prf.\ \ We first notice that if $(z_n)\subset\XX$ is a sequence converging to a point $\tilde\eta\in\rand$, then 
the map
\[ d(x,z_n)-d(\cdot, z_n)\quad\text{converges to the map}\quad \bs_{\tilde\eta}(x, \cdot)\]
uniformly on compact sets. So   
\[ \frac{\e^{-sb(x,z_n)}}{\e^{-sb(y,z_n)}}=\e^{-s(b(x,z_n)-b(y,z_n))}=\e^{-s\delta(\Gamma) (d(x,z_n)-d(y,z_n))}\ \longrightarrow\ \e^{-s\delta(\Gamma)\bs_{\tilde\eta}(x,y)} \]
as $n\to\infty$.  
Hence it suffices to prove that for any sequence $(z_n)\subset\XX$ with $d(\xo,z_n)\to\infty$ we have
\[ \lim_{n\to\infty} \frac{h(\e^{b(x,z_n)})}{h(\e^{b(y,z_n)})}=1.\]
Let $\eps>0$ be arbitrary and fix $\displaystyle\alpha<\frac{\ln(1+\eps)}{\delta(\Gamma) d(x,y)}$. Then by Patterson's Lemma~\ref{Pt} (ii) there exists $r_0>0$  \st  for $r\ge r_0$ and $t>1$
\[  \frac{h(rt)}{h(r)}\le t^{\alpha} .\]
So for all $z\in\XX$ with $b(x,z)\ge \ln(r_0)$ and $b(y,z)\ge \ln(r_0)$ we get
\[   \e^{-\alpha|b(x,z)-b(y,z)|} \le \frac{h(\e^{b(x,z)})}{h(\e^{b(y,z)})}\le \e^{\alpha |b(x,z)-b(y,z)|}.\]
By the remark following the definition of $b$  we have $|b(x,z)-b(y,z)|\le \delta(\Gamma)d(x,y)$ for all $z\in\XX$, hence by choice of $\alpha $
\[ \e^{-\alpha  \delta(\Gamma) d(x,y)} \le \frac{h(\e^{b(x,z)})}{h(\e^{b(y,z)})}\le \e^{\alpha  \delta(\Gamma) d(x,y)}<1+\eps.\]
Using again the inequality $\e^t+\e^{-t}\ge 2$ we obtain as a lower bound
\[ \frac{h(\e^{b(x,z)})}{h(\e^{b(y,z)})} \ge \e^{-\alpha  \delta(\Gamma) d(x,y)} \ge 2- \e^{\alpha  \delta(\Gamma) d(x,y)}>2-(1+\eps)=1-\eps.\]
Hence for all $n\in \NN$ such that $\displaystyle d(\xo,z_n)\ge \frac{\ln(r_0)}{\delta(\Gamma)}+\max\{d(\xo,x),d(\xo,y)\}$ we have\\[2mm]
$\hspace*{5cm} \displaystyle\Big| 1- \frac{h(\e^{b(x,z_n)})}{h(\e^{b(y,z_n)})}\Big|<\eps.$
 \qed\\[1mm]

Recall that $\MM^+(\rand)$ denotes the cone of positive finite Borel measures on $\rand$.
\begin{prp}\label{RadNykPrf}
Let $(s_j)\subset\RR$, $s_j\searrow 1$ be a sequence \st $\mu_\xo^{s_j}$ converges weakly to $\mu_\xo$, and $x\in\XX$ arbitrary.
Then the sequence of measures $\mu_x^{s_j}$ converges weakly to a measure $\mu_x\in{\cal M}^+(\rand)$ with $\supp(\mu_x)\subset\Lim$ and
\[   \frac{d\mu_x}{d\mu_\xo}(\tilde\eta)=\e^{\delta(\Gamma) \bs_{\tilde \eta}(\xo,x)}\qquad\mbox{for all }\ \tilde\eta\in\supp(\mu_{\xo}).\]
\end{prp}
\prf.\ \  Let $f\in \Cnt^0(\ganz)$ with $\Vert f\Vert_\infty<\infty$ be arbitrary,
 $s>1$ and denote 
\mbox{$\,g_{x,\xo}^s:\XX\to\RR$}\break the function defined by (\ref{gxysDef}) in Lemma~\ref{supnorm} (which extends continuously to $\ganz$ by Lemma~\ref{extendscont}).
By~(\ref{RadNikclassical}) we have for all $s>1$ and for all $z\in\supp(\mu_\xo^s)$ 
\[ \frac{\d \mu_x^s}{\d \mu_\xo^s}(z)= \frac{\e^{-sb(x,z)}h(\e^{b(x,z)})}{\e^{-sb(\xo,z)}h(\e^{b(\xo,z)})}=g_{x,\xo}^s(z),\]
hence  
\begin{align*}\int_{\ganz} f(z)\d \mu_x^{s}(z)&= \int_{\ganz} f(z)g_{x,\xo}^s(z)\d \mu_\xo^{s}(z).
\end{align*}
We claim that for any sequence $(s_j)\searrow 1$ \st $\mu_\xo^{s_j}$ converges weakly to $\mu_\xo$ 
we have
\[\tag{$*$} \lim_{j\to\infty} \int_{\ganz} f(z)\d \mu_x^{s_j}(z)=\int_{\ganz} f(\tilde \eta)\e^{\delta(\Gamma)\bs_{\tilde\eta}(\xo,x)} \d \mu_\xo(\tilde\eta);\]
so the measure $\mu_x$ defined by
\[ \frac{\d \mu_x}{\d \mu_\xo}(\tilde \eta)=  \e^{\delta(\Gamma)\bs_{\tilde\eta}(\xo,x)}\qquad\mbox{for all }\ \tilde\eta\in\supp(\mu_{\xo}) \]
is the weak limit of the sequence of measures $\mu_x^{s_j}$. 
Hence in particular we have
\[\supp (\mu_x)\subset\supp (\mu_\xo)\subset\Lim.\]
In order to prove $(*)$ we notice that by Lemma~\ref{extendscont} we have for $\tilde\eta\in\supp(\mu_\xo)\subset\rand$ 
\[ \e^{\delta(\Gamma)\bs_{\tilde\eta}(\xo,x)} = g_{x,\xo}^1(\tilde\eta);\]
we estimate
\begin{align*}
\Big| \int_{\ganz} f(z)\d \mu_x^{s_j}(z) &- \int_{\ganz} f(\tilde \eta)\e^{\delta(\Gamma)\bs_{\tilde\eta}(\xo,x)} \d \mu_\xo(\tilde\eta)\Big| \\
& \le \Big| \int_{\ganz} f(z)g_{x,\xo}^{s_j}(z)\d \mu_\xo^{s_j}(z)
  - \int_{\ganz} f(z)g_{x,\xo}^1(z) \d \mu_\xo^{s_j}(z)\Big|\\
  & + \Big| \int_{\ganz} f(z)g_{x,\xo}^1(z) \d \mu_\xo^{s_j}(z) - \int_{\ganz} f(\tilde \eta) g_{x,\xo}^1(\tilde\eta) \d \mu_\xo(\tilde\eta)\Big|.
\end{align*}
Since $f\cdot g_{x,\xo}^1$ is a bounded and continuous function on $\ganz$, and $\mu_\xo^{s_j}$ converges weakly to $\mu_\xo$, the second term tends to zero as  $j$ tends to infinity. For the
first term we get by definition of the measure $\mu_\xo^{s_j}$
\begin{align*}
\Big| \int_{\ganz} f(z)g_{x,\xo}^{s_j}(z)\d \mu_\xo^{s_j}(z)
 & - \int_{\ganz} f(z)g_{x,\xo}^1(z) \d \mu_\xo^{s_j}(z)\Big|\\
 &=\Big| \int_{\ganz} f(z)\big(g_{x,\xo}^{s_j}(z)- g_{x,\xo}^1(z)\big)\d \mu_\xo^{s_j}(z)\Big|\\
 &=\frac1{\Summ^{s_j}}\sum_{\gamma\in\Gamma} f(\gamma\xo) \big(g_{x,\xo}^{s_j}(\gamma\xo)- g_{x,\xo}^1(\gamma\xo)\big)\e^{-s_j b(\xo,\gamma\xo)}h(\e^{b(\xo,\gamma\xo)})\\
 &\le \Vert f\Vert_\infty \Vert g_{x,\xo}^{s_j}- g_{x,\xo}^1\Vert_\infty \ \longrightarrow\ 0
 \end{align*}
as $j\to\infty$ by Lemma~\ref{supnorm}.\qed\\[1mm]


Recall that $\theta^*\in E$ is the unique unit vector \st $\delta_{\theta^*}(\Gamma)=\delta(\Gamma)$. 
In order to prove Theorem~A it remains to show that the support of the conformal density $\mu$ constructed  above is included in the unique 
$\Gamma$-invariant subset of the limit set which consists of all limit points with slope 
$\theta^*$. For that we need the following auxiliary result which easily follows from
Proposition~\ref{convdiv} (b).
\begin{lem}\label{classconverges}
If $b_\gamma$ is given by~(\ref{bClas}) and $h$ is a non-decreasing function as in Patterson's Lemma~\ref{Pt}, then for all $\eps>0$
$$\sum_{\begin{smallmatrix}{\gamma\in\Gamma}\\
{\Vert\thet(\xo,\gamma\xo)-\theta^*\Vert>\eps}\end{smallmatrix}} \e^{-b_\gamma} h(\e^{b_\gamma}) <\infty.$$
\end{lem}
\prf. \ \ 
Let $\eps >0$ arbitrary and set 
$$s_\eps:=\max\{\delta_\theta(\Gamma): \ \theta\in E, \, \Vert\theta-\theta^*\Vert\ge \eps\}.$$
Then by choice of $\theta^*$ we have $\delta(\Gamma)=\delta_{\theta^*}(\Gamma)>s_\eps$. Fix $\displaystyle \alpha :=\frac12-\frac{s_\eps}{2\delta(\Gamma)}$ and let $r_0>0$ \st for all $r\ge r_0$ 
and $t>1$ we have $h(rt)\le t^{\alpha } h(r)$. In particular, if 
$b_\gamma> \ln(r_0)$, then 
$$h(\e^{b_\gamma})=h\big(\frac{\e^{b_\gamma}}{r_0}\cdot r_0\big)\le \left(\frac{\e^{b_\gamma}}{r_0}\right)^{\alpha }\cdot h(r_0)=\frac{ h(r_0)}{r_0^{\alpha }}\cdot \e^{\alpha  b_\gamma}.$$
Set $\Gamma_\eps:=\{\gamma\in\Gamma:\ \Vert\thet(\xo,\gamma\xo)-\theta^*\Vert>\eps,\, b_\gamma\ge \ln(r_0)\}.$ Then 
\begin{align*}\sum_{\gamma\in\Gamma_\eps}\e^{-b_\gamma} h(\e^{b_\gamma}) &\le  \frac{ h(r_0)}{r_0^{\alpha }} \sum _{\gamma\in\Gamma_\eps}\e^{\alpha b_\gamma}\e^{-b_\gamma}\\
&=  \frac{h(r_0)}{r_0^{\alpha }} \sum _{\gamma\in\Gamma_\eps}\e^{-\delta(\Gamma)(1-\alpha ) d(\xo,\gamma\xo)}.
\end{align*}
 Since $\delta(\Gamma)(1-\alpha ) =\displaystyle\delta(\Gamma)\Big( \frac12+\frac{s_\eps}{2\delta(\Gamma)}\Big)=\frac12\delta(\Gamma)+\frac12 s_\eps >s_\eps$, we conclude that 
$$\sum_{\gamma\in\Gamma_\eps}\e^{-b_\gamma} h(\e^{b_\gamma})\quad\mbox{converges}.$$ The claim now follows from the fact that the set $\{\gamma\in\Gamma:d(\xo,\gamma\xo)\le \ln(r_0)/\delta(\Gamma)\}$ is finite.
\qed\\[1mm]

We finally provide the missing piece in the proof of Theorem A:
\begin{prp}\label{suppconfdens}
The support of the conformal density $\mu=(\mu_x)_{x\in\XX}$ is contained in $\Lim\cap\rand_{\theta^*}$. 
\end{prp}
\prf.\ \ By construction of the map $\mu:\XX\to \MM^+(\rand)$ in Proposition~\ref{RadNykPrf} it suffices  to show that $\supp(\mu_\xo)\subset\Lim\cap\rand_{\theta^*}$.  We already know by definition of $\mu_\xo$ that 
\[\supp(\mu_\xo)\subset\Lim\subset\rand,\]
so it suffices to prove that every point  $\tilde\xi\in \rand\setminus\rand_{\theta^*}$  possesses an open neighborhood $U\subset\ganz$ \st $\mu_\xo(U)=0$. By construction of the measure $\mu_\xo$ as a weak accumulation point of the set
$ \{\mu_\xo^s : s>1\}\subset\MM^+(\ganz)$ for $s\searrow 1$ with $\mu_\xo^s$  defined in~(\ref{probmeasdef}), this is a consequence of Lemma~\ref{classconverges}. \qed

\section{The generalized Patterson-Sullivan construction}\label{genPseries}
According to the statement of Theorem~A, the classical conformal density constructed in the previous section gives measure zero to the set of limit points 
of slope different from $\theta^*$. In order to obtain measures on an arbitrary $\Gamma$-invariant subset of the limit set we will
 use a variation of the classical Patterson-Sullivan construction with more degrees of freedom. The idea is to use a weighted version of the Poincar{\'e} series in order to get the
main contribution from orbit points with direction close to the desired slope $\theta\in E^+$. At this point, properties of the exponent of  growth 
and Proposition~\ref{convdiv} will turn out to be of central importance.

Recall that $\bs_\theta$ denotes the directional distance introduced in Definition~\ref{dirdist}. 
We observe  that for any $b=(b_1, b_2,\ldots,b_r)\in \RR^r$, $\theta\in E$  and
$\tau\ge 0$ fixed, the sum
$$ P_{\theta}^{s,b,\tau}(x,y)=\sum_{\gamma\in \Gamma} \e^{-s
 (b_1 d_1(x_1,\gamma_1 y_1)+b_2 d_2(x_2,\gamma_2 y_2)+\cdots+b_r d_r(x_r,\gamma_r y_r)+\tau(d(x,\gamma y)-\bs_
{\theta}(x,\gamma y )))} $$
possesses an exponent of convergence which is independent
of $x=(x_1,x_2,\ldots,x_r)$, $y=(y_1,y_2,\ldots,y_r)\in\XX$ by the triangle
inequalities  for $d$, $d_1$, $d_2,\ldots, d_r$ and
$\bs_{\theta}$. Notice that for $\tau=0$ this is exactly the sum considered by M.~Burger \cite{MR1230298} in the case of two factors; here we will need to take $\tau>0$  in order to make the 
contribution of  orbit points with direction far away from $\theta$ negligible. 
 
For any  $\theta\in E$ and $\tau\ge 0$, we 
  define a  {\hl region of convergence}
$$ {\mathcal R}_{\theta}^{\tau}:=\big\{b{=}(b_1, b_2,\ldots,b_r) :
\ P_{\theta}^{s,b,\tau}(o,o)\ \text{has exponent of convergence}\
s\le1\}\subset \RR^r$$
and its boundary
$$ \partial {\mathcal R}_{\theta}^{\tau}:=\big\{b{=}(b_1, b_2,\ldots,b_r) :
\ P_{\theta}^{s,b,\tau}(o,o)\ \text{has exponent of convergence}\
s=1\}\subset \RR^r.$$
We recall the definition of the distance vector from (\ref{distancevector}). In the sequel we will identify $b=(b_1,b_2,\ldots,b_r)$ with the column vector $b^\top$ so that for $q=(q_1,q_2,\ldots,q_r)^\top\in\RR^r$ we may write
$$ \langle b,q\rangle =b_1q_1+b_2q_2+\cdots+b_rq_r.$$
The region of convergence  possesses the following two properties:
\begin{lem}
If $\tau\le \tau'$, then ${\mathcal R}_{\theta}^{\tau}\subset {\mathcal R}_{\theta}^{\tau'}$.
\end{lem}
\prf.\ \   Let $\tau\le \tau'$, $b\in {\mathcal
R}_{ \theta}^{\tau}$. Then for any $\gamma\in\Gamma$
$$ \e^{-s\big(\langle b, H(\xo,\gamma\xo)\rangle+\tau'(d(o,\gamma o)-\bs_{ \theta}(o,\gamma
    o))\big)}
{\le}\,\e^{-s\big(\langle b, H(\xo,\gamma\xo)\rangle+\tau(d(o,\gamma o)-\bs_{ \theta}(o,\gamma
    o))\big)} $$
and therefore 
$P_{ \theta}^{s,b,\tau'}(o,o)\le P_{ \theta}^{s,b,\tau}(o,o)$. Hence
  $P_{ \theta}^{s,b,\tau'}(o,o)$   converges if $s>1$. In particular, $P_{ \theta}^{s,b,\tau'}(o,o)$ has
  exponent of convergence less than or equal to 1.\qed

\begin{lem}
For any $\tau\ge 0$, the region ${\mathcal R}_{ \theta}^\tau $ is convex.
\end{lem}
\prf.\ \   Let $\tau\ge 0$, $a, b\in {\mathcal R}_{ \theta}^\tau$ and  $t\in[0,1]$.
For $\gamma\in\Gamma$ we abbreviate 
$$\big(ta+(1-t)b\big)_\gamma:=\langle ta+(1-t)b
, H(\xo,\gamma\xo)\rangle +\tau\big(d(o,\gamma
  o)-\bs_{ \theta}(o,\gamma o)\big).$$
Then by H{\"o}lder's inequality 
$$\sum_{\gamma\in \Gamma} \e^{-s(t a+(1-t)b)_\gamma}=\sum_{\gamma\in
  \Gamma} \e^{-s t a_\gamma} \e^{-s (1-t) b_\gamma}\le \Big(\sum_{\gamma\in
    \Gamma}\e^{-sa_\gamma}\Big)^t  \Big(\sum_{\gamma\in \Gamma}
  \e^{-sb_\gamma}\Big)^{1-t}. $$
The latter sum converges if $s>1$, hence $ta+(1-t)b \in {\mathcal
  R}_{ \theta}^\tau$. \qed\\[-1mm]

With the help of Proposition~\ref{convdiv} we can describe the region of convergence more precisely. The following result 
relates the region of convergence ${\mathcal
  R}^\tau_{ \theta}$ to the exponent of growth of slope $\theta$. 
\begin{lem}\label{bthetgross}
Let $\theta=(\theta_1,\theta_2,\ldots,\theta_r)\in E$   and $\tau\ge 0$. If  $b=(b_1,b_2,\ldots,b_r)\in {\mathcal R}^\tau_{ \theta}$, then
\[  \langle b,\theta\rangle = \sum_{i=1}^r b_i\theta_i \ge \delta_{ \theta}(\Gamma).\]
\end{lem}
\prf.\ \    Assume that 
$ \langle b,\theta\rangle < \delta_{ \theta}(\Gamma)$. Then there exists 
$s>1$ \st $s\langle b,\theta\rangle < \delta_{ \theta}(\Gamma)$. For $H\in\RR^r_{\ge 0}$ we set
$$f(H):=s\big(\langle b, H\rangle+\tau (\Vert H\Vert-\langle H, \theta\rangle)\big),$$
so the continuous homogeneous function $f:\RR_{\ge 0}^r\to\RR$  satisfies  
\[f(\theta)=s\langle b, \theta\rangle  <
  \delta_{ \theta}(\Gamma);\] hence according to  Proposition \ref{convdiv} (a)  applied
to $D=E$, the sum 
$$\sum_{\gamma\in\Gamma} \e^{-f(H(\xo,\gamma\xo))}\qquad\mbox{diverges}.$$
Since 
\vspace{-8mm}
\begin{align*} f(H(\xo,\gamma\xo))=& s\big(b_1 d_1(o_1, \gamma_1 o_1)+b_2 d_2(\xo_2,\gamma_2\xo_2)+\cdots\\ 
&\quad +b_rd_r(\xo_r,\gamma_r\xo_r) +\tau(d(\xo,\gamo)-\bs_{ \theta}(o,\gamo)
)\big)
\end{align*}
we have
\[ \sum_{\gamma\in\Gamma} \e^{-f(H(\xo,\gamma\xo))}=P_\theta^{s,b,\tau}(\xo,\xo),\]
so  we get a contradiction to $(b_1,b_2,\ldots,b_r)\in {\mathcal R}^\tau_{ \theta}$. \qed\\[-1mm]

Using the above properties of the region of convergence and Patterson's Lemma~\ref{Pt} we are now going to construct \bd ies as defined in
the introduction. Such densities are  a natural generalization of
$\Gamma$-invariant conformal densities if one wants to measure an arbitrary   $\Gamma$-invariant subset of the geometric limit set.

From here on we fix $\theta=(\theta_1,\theta_2,\ldots,\theta_r)\in E$ \st $\Lim\cap\rand_\theta\neq \emptyset$, $\tau\ge 0$ and a vector
$b=(b_1,b_2,\ldots,b_r)\in\partial {\mathcal R}^\tau_{ \theta}\subset\RR^r$; let   $\Vert b\Vert{:=}\sqrt{b_1^2+b_2^2+\cdots+b_r^2}\ $ denote its Euclidean norm. 
For $x, z\in \XX$ and 
$\gamma\in\Gamma$   abbreviate  
\begin{align}\label{bgammadef}
b(x,z) &:=  \langle b,H(x,z)\rangle +\tau\big(d(x,z)-\bs_{ \theta}(x,z)\big),\nonumber\\
 b_\gamma &:= b(\xo,\gamma\xo).
\end{align}
Then $\bs_\theta(x,z)\ge 0$ and the  Cauchy--Schwarz inequality 
\[ |\langle b,H(x,z)\rangle |\le \Vert b\Vert\cdot   \Vert H(x,z)\Vert =\Vert b\Vert\cdot d(x,z)\]  
give the important rough estimate
\begin{equation}\label{dbiggerb}
d(x,z)\ge\frac1{\Vert b\Vert+\tau}\cdot b(x,z);
\end{equation}
notice that since the $b_i$ may be negative numbers, 
a converse inequality does not hold in general. In other words, there may exist sequences $(z_n)\subset\XX$ \st $b(\xo,z_n)$ remains bounded even if $d(\xo,z_n)$ tends to infinity. The following lemma gives a condition which ensures that $b(\xo,z_n)$ tends to infinity if $d(\xo,z_n)$ does. 
\begin{lem}\label{bbiggerdifpos} If $\delta_\theta(\Gamma)>0$  and  $(z_n)\subset\XX$ is a sequence converging to a point in $\rand_\theta$, then 
\[b(\xo, z_n)\to\infty\quad\text{as}\quad n\to\infty.\] 
\end{lem}
\prf.\ \ Lemma~\ref{bthetgross} states that for  $b=(b_1,b_2,\ldots,b_r)\in\partial {\mathcal R}^\tau_{ \theta}$ the inequality
\[ \langle b,\theta\rangle \ge \delta_\theta(\Gamma)\] 
holds; since $\delta_\theta(\Gamma)>0$ and the map $\hat\theta\mapsto \langle b, \hat\theta\rangle$ is continuous, there exists $\eps>0$ \st for all $\hat\theta \in E$ with $\Vert\hat\theta-\theta\Vert <\eps$ we have
\[ \langle b,\hat\theta\rangle\ge q>0.\]
Since $(z_n)\subset\XX$ converges to a point in $\rand_\theta$, for all $n$ sufficiently large we have
\[ \Vert\theta-\thet(\xo,z_n)\Vert<\eps,\quad\text{hence}\quad \langle b,H(\xo,z_n)\rangle\ge q\cdot d(\xo,z_n).\]
Summarizing, we get for all $n$ sufficiently large\\[2mm]
$\hspace*{2cm} b(\xo,z_n)=\langle b,H(\xo,z_n)\rangle +\tau\underbrace{\big(d(\xo,z_n)-\bs_{ \theta}(\xo,z_n)\big)}_{\ge 0}\ge q\cdot d(\xo,z_n).$\qed\\[1mm]
Let $h$ be a function as in Patterson's Lemma~\ref{Pt} and recall the definition of the distance vector from~(\ref{distancevector}). 
As in Section~\ref{ClasPatSulconst} we will construct a family of
orbital measures on $\overline{\XX}$ in the following way:
If $D$
denotes the unit Dirac point measure, then  for $x \in\XX$ and $s>1$ we put
\[ \mu_x^s:=\frac1{\Summ^s}\sum_{\gamma\in\Gamma}\e^{-s b(x,\gamma\xo)}h(\e^{b(x,\gamma\xo)})D(\gamma\xo),\quad\text{where}\quad \displaystyle \Summ^s=\sum_{\gamma\in\Gamma}\e^{-sb_\gamma} h(\e^{b_\gamma}).\] 
As in the classical case, these measures are $\Gamma$-equivariant by construction, but now they depend on the additional parameters
$\theta\in E$, $\tau\ge 0$ and $b=(b_1,b_2,\ldots,b_r)\in\partial
{\mathcal R}^\tau_{ \theta}$.  For $x,y\in\XX$ and $s>1$ the measures $\mu_x^s$ and $\mu_y^s$ are absolutely continuous \wrt each other with Radon Nikodym derivative 
\begin{equation}\label{RadNikgeneral} \frac{\d \mu_x^s}{\d \mu_y^s}(z)= \frac{\e^{-sb(x,z)}h(\e^{b(x,z)})}{\e^{-sb(y,z)}h(\e^{b(y,z)})}, \quad z\in\supp(\mu_y^s),
\end{equation}
which again is the function $g_{x,y}^s$ defined by~(\ref{gxysDef}) in Lemma~\ref{supnorm}, but now with the continuous map $b:\XX\times\XX\to\RR\ $ given by~(\ref{bgammadef}). 
For $x=(x_1,x_2,\ldots,x_r),\  y=(y_1,y_2,\ldots,y_r)$,\break $  z=(z_1,z_2,\ldots,z_r)\in\XX$ we further have the estimate
\begin{align}\label{bdiff}
|b(x,z)-b(y,z)| &=\big| \langle b,H(x,z)-H(y,z) \rangle\nonumber \\
&\quad +\tau\big( d(x,z)-d(y,z)\big) -\tau \big( \bs_{ \theta}(x,z)-\bs_{ \theta}(y,z)
\big)\big|\nonumber \\
&\le \Vert b\Vert\cdot \Vert H(x,z)-H(y,z)\Vert +2\tau d(x,y)\le  (\Vert b\Vert +2\tau) d(x,y),
\end{align}
which is~(\ref{bDiffgen}) with constant $C=\Vert b\Vert + 2\tau$. 
So  the conclusion of Lemma~\ref{supnorm} remains true for our new function $g_{x,y}^s:\XX\to\RR$.
Unfortunately an analogous statement of Lemma~\ref{extendscont} does not hold, because in general $g_{x,y}^s$ cannot be extended continuously to the whole geometric boundary. However, the following statement will be sufficient for our purposes: 
\begin{lem}\label{extendscont2} If $\theta\in E^+$ and $\delta_\theta(\Gamma)>0$, then for fixed  $x,y\in \XX$ and $s>0$  the function
\[ g_{x,y}^s:\XX\to\RR,\quad z\mapsto \frac{\e^{-sb(x,z)}h(\e^{b(x,z)})}{\e^{-sb(y,z)}h(\e^{b(y,z)})}\]
extends continuously to  $\XX\cup\rand_\theta$.  Moreover, if 
$x=(x_1,x_2,\ldots,x_r)$, $y=(y_1,y_2,\ldots,y_r)$ and 
$(z_n)\subset\XX$ is a sequence converging to $\tilde\eta=(\eta_1,\eta_2,\ldots,\eta_r,\theta)\in\rand_\theta\subset\regrand$, then
\[ \lim_{n\to\infty} g_{x,y}^s(z_n)=\e^{-s\big(  b_1 \bs_{\eta_1}(x_1,y_1)+b_2 \bs_{\eta_2}(x_2,y_2)+\cdots+b_r\bs_{\eta_r}(x_r,y_r)\big)}.\]
\end{lem}
\prf.\ \ We 
first notice that if $(z_n)=(z_{n,1},z_{n,2},\ldots,z_{n,r})\subset\XX$ is a sequence converging to a point $\tilde\eta=(\eta_1,\eta_2,\ldots,\eta_r,\hat\theta)\in\rand_{\hat\theta}\subset\regrand$ with $\hat\theta=(\hat\theta_1,\hat\theta_2,\ldots,\hat\theta_r)\in E^+$, then  for all $i\in\{1,2,\ldots,r\}$ 
\begin{align*} d_i(x_i,z_{n,i})-d_i(\cdot, z_{n,i})&\quad\text{converges to }\quad \bs_{\eta_i}(x_i, \cdot),\\
\text{and}\quad d(x,z_{n})-d(\cdot, z_{n})&\quad\text{converges to }\quad \bs_{\tilde\eta}(x, \cdot)
\end{align*}
uniformly on compact sets in $\XX_i$ respectively $\XX=\XX_1\times\XX_2\times\cdots\times\XX_r$. Together with Definition~\ref{dirdist} of the directional distance this implies in particular that
\[ \bs_\theta(x,z_n)-\bs_\theta(\cdot ,z_n) \quad\text{converges to }\quad \sum_{i=1}^r \hat\theta_i \cdot \bs_{\eta_i}(x_i, \cdot). \]
Moreover, if $\theta_i\in\RR_{>0}$, $i\in\{1,2,\ldots,r\}$, denotes the $i$-th coordinate of $\theta\in E^+$, then Lemma~\ref{RelProdFactorr} gives
\[ \bs_{\widetilde\eta}(x,y)= \sum_{i=1}^r \theta_i\cdot \bs_{\eta_i}(x_i,y_i).\]
 So  we conclude
\begin{align*}
\lim_{n\to\infty} \big(b(x,z_n)-b(y,z_n)\big) &=\sum_{i=1}^r b_i\bs_{\eta_i}(x_i, y_i)
+\tau\bs_{\widetilde\eta}(x,y)  -\tau\Big(\sum_{i=1}^r \hat\theta_i \cdot\bs_{\eta_i}(x_i,y_i)\Big) \\
&= \sum_{i=1}^r b_i\bs_{\eta_i}(x_i, y_i)
+\tau\Big(\sum_{i=1}^r (\theta_i-\hat\theta_i)\cdot \bs_{\eta_i}(x_i,y_i)\Big),
\end{align*}
and therefore in the case $\hat\theta=\theta$
\begin{align*}
\lim_{n\to\infty} \frac{\e^{-sb(x,z_n)}}{\e^{-sb(y,z_n)}}&=\lim_{n\to\infty} \e^{-s(b(x,z_n)-b(y,z_n))}\\
&=\e^{-s\big(b_1\bs_{\eta_1}(x_1, y_1)+b_2\bs_{\eta_2}(x_2,z_2)+\cdots+b_r\bs_{\eta_r}(x_r,z_r)\big)} .
 \end{align*}
As in the proof of Lemma~\ref{extendscont} we next show that
for any sequence $(z_n)\subset\XX$ with $b(\xo,z_n)\to\infty$ we have
\[ \lim_{n\to\infty} \frac{h(\e^{b(x,z_n)})}{h(\e^{b(y,z_n)})}=1.\]
Let $\eps>0$ be arbitrarily small and fix $\displaystyle\alpha <\frac{\ln(1+\eps)}{(\Vert b\Vert+2\tau) d(x,y)}$. Then by Patterson's Lemma~\ref{Pt} (ii) there exists $r_0>0$  \st  for $r\ge r_0$ and $t>1$
\[  \frac{h(rt)}{h(r)}\le t^{\alpha } .\]
So for all $z\in\XX$ with $b(x,z)\ge\ln( r_0)$ and $b(y,z)\ge\ln( r_0)$ we get
\[   \e^{-\alpha |b(x,z)-b(y,z)|} \le \frac{h(\e^{b(x,z)})}{h(\e^{b(y,z)})}\le \e^{\alpha |b(x,z)-b(y,z)|}.\]
From the estimate~(\ref{bdiff}) and  by choice of $\alpha $
we have 
\[ \e^{- \alpha  (\Vert b\Vert+2\tau)d(x,y)} \le \frac{h(\e^{b(x,z)})}{h(\e^{b(y,z)})}\le \e^{\alpha   (\Vert b\Vert+2\tau) d(x,y)}<1+\eps.\]
Using again the inequality $\e^t+\e^{-t}\ge 2$ we get as a lower bound
\[ \frac{h(\e^{b(x,z)})}{h(\e^{b(y,z)})} \ge \e^{-\alpha   (\Vert b\Vert+2\tau) d(x,y)} \ge 2- \e^{\alpha   (\Vert b\Vert+2\tau) d(x,y)}>2-(1+\eps)=1-\eps.\]
So for $n$ sufficiently large
\[ \Big| 1- \frac{h(\e^{b(x,z_n)})}{h(\e^{b(y,z_n)})}\Big|<\eps;\]
the problem here is that $b(\xo,z_n)$ may remain bounded even if $d(\xo,z_n)$ tends to infinity. However, Lemma~\ref{bbiggerdifpos}  ensures that this  does not happen for sequences  $(z_n)\subset\XX$ converging to a point in $\rand_\theta$. \qed\\[-1mm]

We emphasize again that unlike in the case $b(x,z)=\delta(\Gamma)d(x,z)$ for the classical construction, the continuous function $g_{x,y}^s$ considered here need not extend continuously to the whole geometric boundary. 
One obstruction is the fact that if a sequence\break $(z_n)=(z_{n,1},z_{n,2},\ldots,z_{n,r})$ converges to a point in $\singrand$, then its projections to one or more factors $\XX_i$ need not converge. And even if a sequence $(z_n)$ converges to a regular boundary point  with a slope different from $\theta$,   $b(\xo,z_n)$ may remain bounded and hence the quotient 
$\displaystyle \frac{h(\e^{b(x,z_n)})}{h(\e^{b(y,z_n)})}$ does not necessarily tend to one. 

So in general -- for arbitrary $b=(b_1,b_2,\ldots,b_r)\in\partial
{\mathcal R}^\tau_{ \theta}$ -- there is no analogon of Proposition~\ref{RadNykPrf}. 
However, we still have some freedom in choosing appropriate parameters $b=(b_1,b_2,\ldots,b_r)\in\RR^r$, which can be done as follows: Since the homogeneous extension  $\Psi_\Gamma:\RR_{\ge 0}^r\to \RR$ of the exponent of growth  is concave and upper semi-continuous, it is continuous on the closed convex cone
\[ \ell_\Gamma:=\{ H\in \RR^r_{\ge 0}:\ \Psi_\Gamma(H)\ge 0\} .\] 
So for any $p\in\RR_{\ge 0}^r$ in the relative interior of the intersection of $\ell_\Gamma$ with the vector subspace of $\RR^r$ it spans -- which thanks to Theorem~7.9 in \cite{1308.5584} is equal to the set of points $p\in\RR_{\ge 0}^r$ with  $\Psi_\Gamma(p)>0$ -- 
there exists a linear functional $\Phi$ on $\RR^r$ \st
\begin{equation}\label{tangentfunctional}
\Phi(p)= \Psi_\Gamma(p),\quad\text{ and}\quad\Phi(q)\ge\Psi_\Gamma(q)
\quad\text{for all}\quad q\in\RR_{\ge 0}^r;
\end{equation}
if $\Psi_\Gamma$ is differentiable at $p$, then this linear functional $\Phi$ is unique, but in general it is not. 
For obvious reasons we will call a linear functional $\Phi$ satisfying~(\ref{tangentfunctional}) {\hl tangent} to $\Psi_\Gamma$ at the point $p\in\RR_{\ge 0}^r$. Similarly, we will call a vector $b=(b_1,b_2,\ldots,b_r)\in\RR^r$ {\hl tangent} to $\Psi_\Gamma$ at  $p\in\RR_{\ge 0}^r$, if the  linear functional 
\[ \Phi:\RR^r\to\RR,\quad q\mapsto \langle b, q \rangle\] is tangent to $\Psi_\Gamma$ at the point $p$.  Notice that  if  $b$ is tangent to $\Psi_\Gamma$ at a point $\theta\in E$, then Proposition~\ref{convdiv} implies $b\in \partial
{\mathcal R}^\tau_{ \theta}$. 

It will turn out in the sequel that the choice of $b=(b_1,b_2,\ldots,b_r)\in\RR^r$ tangent to $\Psi_\Gamma$ at $\theta$  is the suitable one.
The following key proposition analogous to Lemma~\ref{classconverges} implies that with this choice any weak accumulation point  of the set of measures $\{\mu_\xo^{s} :s>1\}\subset\MM^+(\ganz)$ as $s\searrow 1$ is supported in $\rand_\theta$.  It is therefore the key ingredient 
in 
the construction of orbital measures with support in a single $\Gamma$-invariant subset $ \Lim\cap \rand_\theta\subset\rand$. 
\begin{prp}\label{key}
Fix $\theta\in E^+$ \st $\delta_\theta(\Gamma)>0$,  and let   $b=(b_1,b_2,\ldots,b_r)\in\RR^r$ be a vector tangent to $\Psi_\Gamma$ at $\theta$. Then for  all $\tau>0$ and for 
all $ \eps>0$
\[\sum_{\begin{smallmatrix}{\gamma\in\Gamma}\\{\Vert\thet(\xo,\gamma\xo)- \theta\Vert> \eps}\end{smallmatrix}}
\e^{-b_\gamma} h(\e^{b_\gamma}) <\infty,\]
where $b_\gamma$ is defined in~(\ref{bgammadef}) and $h$ is a function as in Patterson's Lemma~\ref{Pt}.
\end{prp}
\prf.\ \ Since $b=(b_1,b_2,\ldots,b_r)\in\RR^r$ is tangent to $\Psi_\Gamma$ at $\theta$ and
$\Psi_\Gamma(\theta)= \delta_\theta(\Gamma)$ we have
$$ \langle b,\theta\rangle =\delta_\theta(\Gamma)\quad\text{ and}\quad \langle b,\hat\theta\rangle\ge \delta_{\hat\theta}(\Gamma)\quad\text{ for all}\quad \hat\theta\in E.$$ 
We fix $\tau>0$ and let $\eps >0$ be arbitrary. Since the sum is non-increasing when $\eps$ gets bigger, we may further assume that $\eps <1$. 

We set $D:=\{\hat\theta\in E:\ \Vert \hat\theta-\theta\Vert >\eps\}$,
$s_\eps:=\max\{|\langle b,\hat\theta\rangle|  : \ \hat \theta\in \overline{D}\}\ge 0 $ ,
and fix  
\[\alpha <\frac{\tau\eps^2}{\tau\eps^2+2 s_\eps} \le 1.\]
By Patterson's Lemma~\ref{Pt} (ii) there exists $r_0>0$ \st for all $t>1$ we have
$h(r_0t)\le t^{\alpha } h(r_0)$. So $b_\gamma> \ln(r_0)$ implies
\[ h(\e^{b_\gamma})=h\Big(\frac{\e^{b_\gamma}}{r_0}\cdot r_0\Big)\le\Big( \frac{\e^{b_\gamma}}{r_0}\Big)^{\alpha } h(r_0).\]
We therefore have
\begin{align*}
\sum_{\begin{smallmatrix}{\gamma\in\Gamma}\\{\Vert\thet(\xo,\gamma\xo)- \theta\Vert> \eps}\\{b_\gamma> \ln(r_0)}\end{smallmatrix}}
\e^{-b_\gamma} h(\e^{b_\gamma}) &\le \frac{h(r_0)}{r_0^{\alpha }} \sum_{\begin{smallmatrix}{\gamma\in\Gamma}\\{\Vert \thet(\xo,\gamma\xo)- \theta\Vert> \eps}\\{b_\gamma> \ln(r_0)}\end{smallmatrix}}
\e^{-(1-\alpha )b_\gamma} <\infty,
\end{align*}
by Proposition~\ref{convdiv}. Indeed,   the continuous homogeneous function $f:\RR_{\ge 0}^r\to\RR$ defined by
\[ f(H):= (1-\alpha )\big(\langle b,H\rangle +\tau(\Vert H\Vert -\langle H,\theta\rangle)\big)\]
satisfies  $f(H(\xo,\gamma\xo))=(1-\alpha )b_\gamma\ $ for all $\gamma\in\Gamma$, and 
\[ f(\hat\theta)>\delta_{\hat\theta}(\Gamma)\quad\text{for all}\quad  \hat\theta\in\overline{D}=\{\hat\theta\in E:\ \Vert \hat\theta-\theta\Vert \ge \eps\}\] by the following estimate:
\begin{align*}
f(\hat\theta)&=(1-\alpha )\big(\langle b,\hat\theta\rangle +\tau(1 -\langle\hat\theta,\theta\rangle)\big) =(1-\alpha )\langle b,\hat\theta\rangle +(1-\alpha )\frac{\tau}2\underbrace{\Vert \hat\theta-\theta\Vert^2}_{\ge \eps^2}\\
&\ge \langle b,\hat\theta\rangle -\alpha  |\langle b,\hat\theta\rangle|+(1-\alpha )\tau\frac{\eps^2}{2}\ge  \langle b,\hat\theta\rangle +\tau\frac{\eps^2}{2} - \alpha \big(  |\langle b,\hat\theta\rangle| + \tau\frac{\eps^2}{2} \big)\\
&> \langle b,\hat\theta\rangle +\tau\frac{\eps^2}{2} - \frac{\tau\eps^2}{\tau\eps^2+2 s_\eps} \big(  \underbrace{|\langle b,\hat\theta\rangle|}_{\le s_\eps} + \tau\frac{\eps^2}{2} \big)\ge \langle b,\hat\theta\rangle\ge \delta_{\hat\theta}(\Gamma).
\end{align*}
So it remains to show that the sum
\[\sum_{\begin{smallmatrix}{\gamma\in\Gamma}\\{\Vert\thet(\xo,\gamma\xo)- \theta\Vert> \eps}\\{b_\gamma \le \ln(r_0)}\end{smallmatrix}}
\e^{-b_\gamma} h(\e^{b_\gamma})\] 
is finite. This is not as trivial as in Lemma~\ref{classconverges}, because due to the fact that the coordinates $b_i$, $i\in\{1,2,\ldots,r\}$, need not be positive,
$ b_{\gamma_n} $  might remain bounded even if $d(\xo,\gamma_n\xo)$ tends to infinity.  However, we can argue in the following way: By 
Property~(a) of the exponent of growth the set of slopes of limit points 
\[P_\Gamma:=\{\theta'\in E :\ \Lim\cap\rand_{\theta'}\ne\emptyset\}\]
satisfies 
\[ P_\Gamma=\{\theta'\in E :\ \delta_\theta(\Gamma)\ge 0\}. \]
Since the map $E\to\RR$, $\hat\theta\mapsto \langle b,\hat \theta\rangle$ is continuous and satisfies $\langle b,\hat\theta\rangle \ge \delta_{\hat\theta}(\Gamma)$ for all $\hat\theta\in E$, there exists  $\alpha >0$ \st for all $\gamma\in\Gamma$ with $\thet(\xo,\gamma\xo)$ in the $\alpha$-neighborhood of $P_\Gamma$ defined by
\[ \{\hat\theta\in E:\ \Vert\hat \theta-\theta'\Vert<\alpha\quad\text{ for some}\quad \theta'\in P_\Gamma\}\] we have
\[\langle b, \thet(\xo,\gamma\xo)\rangle >-\frac{\tau\eps^2}4.\]
Hence for all $\gamma\in\Gamma$ with $\Vert\thet(\xo,\gamma\xo)- \theta\Vert> \eps$ and 
$\Vert\thet(\xo,\gamma\xo)-\theta'\Vert<\alpha $ for some $\theta'\in P_\Gamma$ we get
\begin{align*} b_\gamma &= d(\xo,\gamma\xo)\big(\langle b,\thet(\xo,\gamma\xo)\rangle +\tau\underbrace{(1 -\langle\thet(\xo,\gamma\xo),\theta\rangle)}_{\ge \eps^2/2}\big)\\
&> d(\xo,\gamma\xo)\big(-\frac{\tau\eps^2}4 +\tau\frac{\eps^2}2\big)=\frac{\tau\eps^2}4d(\xo,\gamma\xo).
\end{align*}
So $b_\gamma\le \ln(r_0)$ implies $d(\xo,\gamma\xo)\le\frac4{\tau\eps^2} \ln(r_0)$, hence the
set 
\[ \{\gamma\in\Gamma : \Vert \thet(\xo,\gamma\xo)- \theta\Vert> \eps,\  \Vert\thet(\xo,\gamma\xo)-\theta'\Vert<\alpha  \ \text{for some}\ \theta'\in P_\Gamma, \ b_\gamma \le  \ln(r_0)\}\]
is finite. 

Now assume the remaining set
\[ \Gamma':= \{\gamma\in\Gamma : \Vert\thet(\xo,\gamma\xo)- \theta\Vert> \eps,\  \Vert\thet(\xo,\gamma\xo)-\theta'\Vert\ge \alpha  \ \text{for all}\ \theta'\in P_\Gamma, \ b_\gamma \le \ln(r_0)\}\]
is infinite. Then there exists an accumulation point of the orbit $\Gamma'\xo$ which cannot belong to the geometric boundary $\rand$ of $\XX$ by definition of 
the set $P_\Gamma$. Since $\Gamma'$ is discrete, the accumulation point cannot belong to $\XX$ either, so the only possibility is that  $\Gamma'$ is finite. \qed\\[1mm]
 
 Notice that the proof above also works in the case $\theta=(\theta_1,\theta_2,\ldots,\theta_r)\in E\setminus E^+$. However, in general the vector $b=(b_1,b_2,\ldots,b_r)$ tangent to $\Psi_\Gamma$ at $\theta$ does not satisfy $b_i=0$ for all $i\in\{1,2,\ldots,r\}$ with $\theta_i=0$, which would be necessary for a well-defined Radon Nikodym derivative of the $(b,\theta)$-density. 
Hence for the construction of \bd ies we have to restrict ourselves to $\theta\in E^+$.

From here on we therefore fix a slope $\theta=(\theta_1,\theta_2,\ldots,\theta_r)\in E^+$ such that\break $\Psi_\Gamma(\theta)=\delta_\theta(\Gamma)>0$; hence in particular we have $\Lim\cap\rand_\theta\neq\emptyset$.   
Using the previous proposition we are finally able to construct a \bd y according to Definition~\ref{bdensi}.

To that end we further choose  
$b=(b_1,b_2,\ldots,b_r)\in\RR^r$ tangent to $\Psi_\Gamma$ at $\theta$ and fix \mbox{$\tau>0$.}  
By Proposition~\ref{key} and arguments analogous to the proof of Proposition~\ref{suppconfdens} any weak accumulation point $\mu_\xo$ of  $\mu_\xo^s$ (defined with these parameters) as  $s\searrow 1$  satisfies $\supp(\mu_\xo)\subset \Lim\cap\rand_\theta$. 

The last step in the proof of Theorem~B from the introduction needs an analogon of Proposition~\ref{RadNykPrf}; we have to show that a suitable choice of a weak accumulation point  for each of the sets
\[ \{\mu_x^s:s\in (1,2]\}\subset\MM^+(\ganz),\quad x\in\XX,\]
produces a \bd y.
We recall that for $\tilde\eta\in\rand$ the Busemann vector $\bv_{\tilde\eta}$ was defined in~(\ref{busvector}). 

\begin{prp}\label{RadNykPrfgeneral}
Let $(s_j)\subset\RR$, $s_j\searrow 1$ be a sequence \st $\mu_\xo^{s_j}$ converges weakly to $\mu_\xo$, and $x\in\XX$ arbitrary.
Then the sequence of measures $\mu_x^{s_j}$ converges weakly to a measure $\mu_x\in{\cal M}^+(\rand)$ with $\supp(\mu_x)\subset\Lim\cap\rand_\theta$ and
\[   \frac{d\mu_x}{d\mu_\xo}(\tilde\eta)=\e^{\langle b, \bv_{\tilde \eta}(\xo,x)\rangle}\qquad\mbox{for all }\ \tilde\eta\in\supp(\mu_{\xo}).\]
\end{prp}
\prf.\ \  Let $f\in \Cnt^0(\ganz)$ with $\Vert f\Vert_\infty<\infty$ be arbitrary,
 $s>1$ and denote 
\mbox{$\,g_{x,\xo}^s:\XX\to\RR\,$} the continuous function defined by (\ref{gxysDef}) in Lemma~\ref{supnorm}, which extends continuously to $\XX\cup\rand_\theta$ by Lemma~\ref{extendscont2}.
By~(\ref{RadNikgeneral}) we have for all $s>1$ and for all $z\in\supp(\mu_\xo^s)$ 
\[\frac{\d \mu_x^s}{\d \mu_y^s}(z)= \frac{\e^{-sb(x,z)}h(\e^{b(x,z)})}{\e^{-sb(y,z)}h(\e^{b(y,z)})}=g_{x,\xo}^s(z),\]
hence  
\begin{align*}\int_{\XX\cup\rand_\theta} f(z)\d \mu_x^{s}(z)&= \int_{\XX\cup\rand_\theta} f(z)g_{x,\xo}^s(z)\d \mu_\xo^{s}(z).
\end{align*}
We claim that for any sequence $(s_j)\searrow 1$ \st $\mu_\xo^{s_j}$ converges weakly to $\mu_\xo$ 
we have
\[\tag{$*$} \lim_{j\to\infty} \int_{\ganz} f(z)\d \mu_x^{s_j}(z)=\int_{\ganz} f(\tilde \eta)\e^{\langle b, \bv_{\tilde\eta}(\xo,x)\rangle} \d \mu_\xo(\tilde\eta);\]
so the measure $\mu_x$ defined by
\[ \frac{\d \mu_x}{\d \mu_\xo}(\tilde \eta)=  \e^{\langle b,\bv_{\tilde\eta}(\xo,x)\rangle}\qquad\mbox{for all }\ \tilde\eta\in\supp(\mu_{\xo}) \]
is the weak limit of the sequence of measures $\mu_x^{s_j}$. Then in particular we have
\[\supp (\mu_x)\subset\supp (\mu_\xo)\subset\Lim\cap\rand_\theta.\]
In order to prove $(*)$ we notice that by Lemma~\ref{extendscont2} we have for $\tilde\eta\in\supp(\mu_\xo)\subset\rand_\theta$ 
\[ \e^{\langle b,\bv_{\tilde\eta}(\xo,x)\rangle} = g_{x,\xo}^1(\tilde\eta);\]
moreover, using the fact that the support of any measure in the weak closure of the set 
$\{\mu_x^s: s\in (1,2]\}\subset\MM^+(\ganz)$ is contained in $\XX\cup \rand_\theta$,  we estimate
\begin{align*}
\Big| \int_{\ganz} f(z)\d \mu_x^{s_j}(z) &- \int_{\ganz} f(\tilde \eta)\e^{\langle b,\bv_{\tilde\eta}(\xo,x)\rangle} \d \mu_\xo(\tilde\eta)\Big| \\
&=\Big| \int_{\XX\cup\rand_\theta} f(z)\d \mu_x^{s_j}(z) - \int_{\rand_\theta} f(\tilde \eta)\e^{\langle b,\bv_{\tilde\eta}(\xo,x)\rangle} \d \mu_\xo(\tilde\eta)\Big| \\
& =\Big| \int_{\XX\cup\rand_\theta} f(z)g_{x,\xo}^{s_j}(z)\d \mu_\xo^{s_j}(z) - \int_{\rand_\theta} f(\tilde \eta)g_{x,\xo}^1(\tilde\eta) \d \mu_\xo(\tilde\eta)\Big|\\
& \le \Big| \int_{\XX\cup\rand_\theta} f(z)g_{x,\xo}^{s_j}(z)\d \mu_\xo^{s_j}(z)
  - \int_{\XX\cup\rand_\theta} f(z)g_{x,\xo}^1(z) \d \mu_\xo^{s_j}(z)\Big|\\
  & + \Big| \int_{\XX\cup\rand_\theta} f(z)g_{x,\xo}^1(z) \d \mu_\xo^{s_j}(z) - \int_{\XX\cup\rand_\theta} f(\tilde \eta) g_{x,\xo}^1(\tilde\eta) \d \mu_\xo(\tilde\eta)\Big|.
\end{align*}

Since $f\cdot g_{x,\xo}^1$ is a bounded and continuous function on $\XX\cup\rand_\theta$, and $\mu_\xo^{s_j}$ converges weakly to $\mu_\xo$, the second term tends to zero as  $j$ tends to infinity.

For the
first term we argue as in the proof of Proposition~\ref{RadNykPrf}. By the estimate~(\ref{bdiff}) the assumption~(\ref{bDiffgen}) on $b$ in Lemma~\ref{supnorm} is satisfied with $C=\Vert b\Vert+2\tau$, hence  by definition of the measure $\mu_\xo^{s_j}$ and by Lemma~\ref{supnorm} we have
\begin{align*}
\Big| \int_{\XX\cup\rand_\theta} f(z)g_{x,\xo}^{s_j}(z)\d \mu_\xo^{s_j}(z)
 & - \int_{\XX\cup\rand_\theta} f(z)g_{x,\xo}^1(z) \d \mu_\xo^{s_j}(z)\Big|\\
 &=\Big| \int_{\XX\cup\rand_\theta} f(z)\big(g_{x,\xo}^{s_j}(z)- g_{x,\xo}^1(z)\big)\d \mu_\xo^{s_j}(z)\Big|\\
 &=\frac1{\Summ^{s_j}}\sum_{\gamma\in\Gamma} f(\gamma\xo) \big(g_{x,\xo}^{s_j}(\gamma\xo)- g_{x,\xo}^1(\gamma\xo)\big)\e^{-s_j b(\xo,\gamma\xo)}h(\e^{b(\xo,\gamma\xo)})\\
 &\le \Vert f\Vert_\infty \Vert g_{x,\xo}^{s_j}- g_{x,\xo}^1\Vert_\infty \ \longrightarrow\ 0\quad\text{as}\quad j\to\infty.
 \end{align*}\qed\\[1mm]


So given a slope $\theta\in E^+$ satisfying $\Psi_\Gamma(\theta)=\delta_\theta(\Gamma)>0$ our construction with $b=(b_1,b_2,\ldots,b_r)\in\RR^r$ tangent to $\Psi_\Gamma$ at the point $\theta$ and $\tau>0$ arbitrary produces the desired $(b,\theta)$-density.  
 This proves Theorem~B
from the introduction.

\section{Properties of $(b,\theta)$-densities}\label{Propbdies}

In this section we will study properties of \bd ies using the shadow lemma Theorem~\ref{shadowlemma}. If not otherwise specified we allow $\theta\in E$.

\begin{lem}
Let $\mu$ be a \bd y, and $x\in\XX$.  If $\,\widetilde U\subset\rand$ is an open neighborhood of a limit point $\tilde\xi\in\rand_\theta$, then  
$\mu_x(\widetilde U) >0$. 
\end{lem}
\prf. \ \ Let $\widetilde U \subset\rand$ be an  open neighborhood of a limit point $\tilde\xi\in\rand_\theta$ \st $\mu_x(\widetilde U)=0$. If $U:=\widetilde U\cap\rand_\theta$, then by compactness 
and minimality of $\Lim\cap\rand_\theta$ (see Theorem~A in \cite{1308.5584}) 
there exists a finite set $\Lambda\subset\Gamma$ \st 
$$\Lim\cap\rand_\theta  \subset \bigcup_{\gamma\in\Lambda} \gamma U.$$
Moreover, by $\Gamma$-equivariance
$$ \mu_x(\Lim\cap\rand_\theta)\le\sum_{\gamma\in\Lambda}\mu_x(\gamma U)=\sum_{\gamma\in\Lambda}\mu_{\gamma^{-1}x}(U)\le \sum_{\gamma\in\Lambda}\mu_{\gamma^{-1}x}(\widetilde U)=0,$$
since $\mu_{\gamma^{-1}x}$, $\gamma\in\Lambda$, is absolutely continuous with respect to $\mu_x$.\qed\\[-1mm]

Recall the definition of the distance vector (\ref{distancevector}) from Section~\ref{prodHadspaces}. 
\begin{thr}[Shadow lemma]\label{shadowlemma} 
Let $\mu$ be a \bd y. Then there exists a constant $c_0>0$ \st
for any $c>c_0$ 
there exists a constant  $D(c)>1$ with the property 
$$ \frac1{D(c)}\e^{-\langle b, H(\xo,\gamma\xo)\rangle} \le \mu_{o
}\big(\Sh(o:B_{\gamma o}(c))\big)\le D(c)\e^{-\langle b, H(\xo,\gamma\xo)\rangle}$$
for all $\gamma\in\Gamma$ with $d(\xo,\gamma\xo)>c$.
\end{thr}   
\prf. \ \   For $i\in\{1,2,\ldots,r\}$ we let $U_i\subset\rand_i$ be open neighborhoods of $h_i^+$,  $\Lambda\subset\Gamma$ a finite set and $c_0>0$ such that the assertion of 
Proposition~\ref{largeshadows} holds. If $U_\theta$ denotes the Cartesian product of the sets $U_i$ with $i\in I^+(\theta)$, then for all $\lambda\in\Lambda$ the  set 
\[ \lambda(U_\theta\times\{\theta\}) \subset\rand_\theta\]
is a relatively open neighborhood of a limit point in $\rand_\theta$, so by the previous lemma   
\[ q:=\min\{ \mu_\xo\big(\lambda(U_\theta\times\{\theta\})\big): \lambda\in\Lambda\}\]
is strictly positive. Moreover, if $c\ge c_0$ and $\gamma\in\Gamma$ \st $d(o,\gamo)>c$ then by Proposition~\ref{largeshadows} there exists $\lambda\in\Lambda$ \st $\lambda(U_\theta\times\{\theta\})\subset \Sh(\gamma^{-1}\xo:B_\xo(c))$. Hence for $c\ge c_0$ and
$\gamma\in\Gamma$ with $d(o,\gamo)>c$ we have
\begin{equation}\label{mush}
\mu_{o}(\rand)\ge\mu_{o}\big(\Sh(\gamma^{-1}o :B_{o}(c))
\big)\ge q>0.\end{equation}
Put $S_\gamma:=\Sh(\xo:B_{\gamma o}(c))$ and recall the definition of the Busemann vector (\ref{busvector}). The properties (ii) and (iii) of a \bd y imply
\begin{align*}
\mu_{o}\big(\Sh(\gamma^{-1}\xo : B_{o}(c))\big)&=\mu_{o}(\gamma^{-1}
S_\gamma) =\mu_\gamo(S_\gamma)\\
&= \int_{S_\gamma}d\mu_{\gamo}(\tilde\eta)=\int_{\smSh(o:B_{\gamo}(c))}\e^{\langle b,\bv_{\tilde\eta}(\xo,\gamma\xo)\rangle} d\mu_{o}(\tilde\eta).
\end{align*}
By {\rm Lemma}~\ref{esti},
$$ 
\e^{-2c}\e^{\langle b, H(\xo,\gamma\xo)\rangle }\mu_{o}(S_\gamma)
< \mu_{o}\big(\Sh(\gamma^{-1}o:B_{o}(c))\big)
\le  \e^{\langle b, H(\xo,\gamma\xo)\rangle} \mu_{o}(S_\gamma),
$$
so equation (\ref{mush}) allows us to conclude
$$
\e^{-\langle b, H(\xo,\gamma\xo)\rangle} q
\le \mu_{o}(S_\gamma)\le \e^{-\langle b, H(\xo,\gamma\xo)\rangle} \e^{2c}\at\mu_{o}(\rand).\eqno{\scriptstyle\square}
$$

The following applications of Theorem~\ref{shadowlemma} yield relations between the  exponent of growth of a given slope 
$\theta\in E$ and the parameters of a \bd y. 

\begin{thr}\label{deltaklein}
If for $\theta\in E^+$ a $\,\Gamma$-invariant
$(b,\theta)$-density exists, then
$$\delta_{\theta} (\Gamma)\le \langle b, \theta\rangle.$$
\end{thr}
\prf. \ \   Suppose $\mu$ is a \bd y. Let
$c>c_0+1$, where $c_0>0$ is as in {\rm Theorem}~\ref{shadowlemma}, $\eps>0$ and $n\in\NN$, $n>3c_0$  arbitrary. Let $\tilde\eta=(\eta_1,\eta_2,\ldots,\eta_r,\theta)\in\supp(\mu_\xo)$.
We only need $N(\eps) n^{r-1}$ balls of radius $1$ in
$\XX$ to cover the set 
\[\{\big(\sigma_{\xo_1,\eta_1}(t \hat\theta_1),\sigma_{\xo_2,\eta_2}(t\hat\theta_2),\ldots,\sigma_{\xo_r,\eta_r}(t\hat\theta_r)\big)\in\XX: \ n-1\le t <n,\, \Vert\hat\theta -\theta\Vert<\eps\},\]
and
$N(\eps)$ is independent of $n$. Since $\Gamma$ is discrete, a
$2c$-neighborhood of any of these balls contains a uniformly bounded number
$M_c$ of ele\-ments of $\Gamma\cdot\xo$. Hence 
every point in the support of $\mu_\xo$ is contained
%
in at most $M_c N(\eps) n^{r-1} $ Weyl chamber shadows\break
$\Sh(\xo: B_\gamo(c))$ with  $\gamma\in \Gamma':=\{\gamma\in\Gamma :\Vert \thet(\xo,\gamma\xo)-\theta\Vert<\eps,\,n-1\le
    d(\xo,\gamma\xo)<n\}$.
Therefore 
\begin{align*}
\sum_{\gamma\in\Gamma'}\mu_\xo\big(\Sh(\xo: B_\gamo(c))\big) &
\le  M_c  N(\eps) n^{r-1}\cdot \mu_\xo\big(\bigcup_{\gamma\in\Gamma'}\Sh(\xo: B_\gamo(c)) \big)\\
&\le  M_c     N(\eps) n^{r-1} \cdot\mu_\xo(\rand_\theta)
= M_c  N(\eps) n^{r-1} \cdot\mu_\xo(\rand).
\end{align*}
Furthermore, if $\, \gamma\in\Gamma'$ then $\thet(\xo,\gamma\xo)\in E$ 
satisfies $ \Vert \thet(\xo,\gamma\xo) -\theta\Vert\le\eps$. 
Using the Cauchy--Schwarz inequality 
we get for $\gamma\in\Gamma'$
$$
\langle b,\thet(\xo,\gamma\xo)\rangle =\langle b,\theta\rangle +\langle b, \thet(\xo,\gamma\xo)-\theta\rangle\le \langle b,\theta\rangle   +\Vert b\Vert \eps.$$
With the notation 
$$\Delta N_\theta^\eps(\xo,\xo;n):=\#\{ \gamma\in\Gamma:\ n-1< d(\xo,\gamma \xo)\le n,\ \Vert\thet(\xo,\gamma\xo)-\theta\Vert<\eps\} $$ 
from Section~\ref{ExpGrowth} and the shadow lemma  {\rm Theorem}~\ref{shadowlemma} we conclude
\begin{multline*}
\Delta N_{\theta}^\eps(o,o;n)\frac1{D(c)}\e^{-\langle b,\theta\rangle n}{\le}\sum_{\gamma
\in\Gamma'}\frac1{D(c)}
    \e^{-\langle b, H(\xo,\gamma\xo)\rangle +\eps \Vert b\Vert  d(\xo,\gamo)}\\
 \le  \e^{\eps \Vert b\Vert
    n}\sum_{\gamma\in\Gamma'}\mu_\xo\big(\Sh(\xo:B_\gamo(c))\big) 
\le \e^{\eps \Vert b\Vert
    n}  M_c 
    N(\eps) n^{r-1}\cdot \mu_\xo(\rand).
\end{multline*}
Using~(\ref{defbylimsup}) we therefore get 
\begin{align*}
\delta_{\theta}^\eps(o,o)& \le  \limsup_{n\to\infty}
\frac1{n}\ln \big( D(c)  M_c 
N(\eps)  n^{r-1}\cdot\mu_\xo(\rand) \e^{\langle b,\theta\rangle n+\eps \Vert b\Vert
    n}\big)= \langle b,\theta\rangle + \eps \Vert b\Vert
\end{align*}
and the claim follows as $\eps\searrow 0$.\qed\\[-1mm]

We remark that  the proof of the above proposition does not work for $\theta\in E\setminus E^+$ because a singular boundary point can be contained in infinitely many Weyl chamber shadows 
$\Sh(\xo: B_\gamo(c))$ with  $\gamma\in \{\gamma\in\Gamma :\Vert\thet(\xo,\gamma\xo)-\theta\Vert<\eps,\,n-1\le
    d(\xo,\gamma\xo)<n\}$.  This is due to the fact that the Weyl chambers ${\cal C}_{\xo,\tilde\eta}$ with $\tilde\eta\in\singrand$ are too big. 

We next recall the  notion of radial limit point from Definition~\ref{raddef} of the introduction. If $\theta\in E$  then using the equality~(\ref{shadowsingular}) we can describe the radial limit set in $\rand_\theta$ via
\begin{equation}\label{radlimthetacover}
 \radlim\cap\rand_\theta=\bigcup_{c>0}\bigcap_{R>c}\, \bigcap_{\eps>0}
\bigcup_{\begin{smallmatrix}{\gamma\in
\Gamma}\\
{ d(o,\gamma
\xo)>R}\\{ \Vert\thet(\xo,\gamma\xo)-\theta\Vert<
\eps}\end{smallmatrix}}
\Sh(\xo: B_{\gamma\xo}(c))\cap\rand_\theta  
\end{equation}

Together with the previous theorem the following statement says that if a\break \mbox{\bd y} gives positive measure to
the regular radial limit set, then the exponent of growth of $\Gamma$ of slope $\theta$ is completely determined by the parameters $b=(b_1,b_2,\ldots, b_r)$.
\begin{thr}\label{posmeasrad}
If $\theta\in E$ and a \bd y gives positive measure to $\radlim$, then 
$ \delta_{\theta}(\Gamma) \ge  \langle b, \theta\rangle$.
\end{thr}
\prf. \ \   Suppose  $\mu$ is a \bd y  \st $\mu_\xo(\radlim)>0$. Since $\supp(\mu_\xo)\subset\rand_\theta$ this implies $\mu_\xo(\radlim\cap\rand_\theta)>0$.
By definition~(\ref{radlimthetacover}) there exists $c>0$ \st with
\begin{equation}\label{radlimkleinerc}
 \radlim (c)\cap\rand_\theta:=\bigcap_{R>c}\,\bigcap_{\eps>0}\,\bigcup_{\begin{smallmatrix}{\gamma\in
\Gamma}\\
{ d(o,\gamma
\xo)>R}\\{ \Vert\thet(\xo,\gamma\xo)-\theta\Vert<
\eps} \end{smallmatrix}}
\Sh(\xo: B_{\gamma\xo}(c))  \cap\rand_\theta
\end{equation}
we have $\mu_\xo(\radlim(c)\cap\rand_\theta)>0$. 
Without loss of generality we may assume that $c>c_0$  with $c_0>0$ as in Theorem~\ref{shadowlemma}.
Let $\eps>0$ and $R>c$ arbitrary,  and set 
\[\Gamma':=\{\gamma\in\Gamma :\  d(o,\gamma\xo)>R,\ \Vert \thet(\xo,\gamma\xo)-\theta\Vert<\eps\}.\] 
Then by (\ref{radlimkleinerc})
$$\radlim (c)\cap\rand_\theta \subset \bigcup_{\gamma\in\Gamma'}
\Sh(\xo: B_\gamo (c))\cap\rand_\theta,$$
and  we estimate 
\begin{align*}
 0 <\mu_{\xo}(\radlim (c))&=\mu_\xo(\radlim (c)\cap \rand_\theta)\\
 &\le\sum_{\gamma\in\Gamma'}\mu_\xo\big(\Sh(\xo: B_\gamo(c))\big) \le D(c)\sum_{\gamma\in\Gamma'}\e^{-\langle b, H(\xo,\gamma\xo)\rangle}.
 \end{align*}
This  implies that for any
$\eps>0$ the tail of the series
$$
\sum_{\begin{smallmatrix}{\scriptstyle \gamma\in
\Gamma}\\{\scriptstyle \Vert\thet(\xo,\gamo)-\theta\Vert<\eps}\end
{smallmatrix}}\e^{-\langle b, H(\xo,\gamma\xo)\rangle}$$ 
does not tend to zero. Therefore the sum above diverges,
and by {\rm Proposition} \ref{convdiv} (b) there exists
$\hat\theta\in E$, $\Vert\hat\theta-\theta\Vert\le\eps$ \st
$$\langle b,\hat\theta\rangle \le
\delta_{\hat\theta}(\Gamma).$$ 
Taking the limit as $\eps\searrow 0$, we conclude 
$\quad \langle b,\theta\rangle \le\delta_{\theta}
(\Gamma)$. \qed\\[-1mm]

Recall the definition of the Busemann vector (\ref{busvector}) from Section~\ref{prodHadspaces}. The following two lemmata hold for any 
 $\theta\in E$ 
and will be important for the proof of Theorem~\ref{radpoint}.
\begin{lem}\label{busnull}
Let $\mu$ be a $( b,\theta)$-density. If $\tilde\eta\in \rand_\theta$ is a point mass for $\mu$,  and 
$\Gamma_{\tilde\eta}$ its stabilizer, then for 
any $\gamma\in\Gamma_{\tilde\eta}$ and $x=(x_1,x_2,\ldots,x_r)\in \XX\;$  we
have 
\[ \langle b,\bv_{\tilde\eta}(x,\gamma x)\rangle = \sum_{i=1}^r b_i\bs_{\eta_i}(x_i,\gamma_i x_i)=0.\]
In particular, if $\gamma,\hat\gamma\in\Gamma$ are representatives of the same
coset in $\Gamma/\Gamma_{\tilde\eta}$, then 
 $$\langle b,\bv_{\tilde\eta}(x,\gamma^{-1} x)\rangle= \langle b,\bv_{\tilde\eta}(x,\hat\gamma^{-1} x)\rangle .$$
\end{lem}
\prf. \ \  For $x\in\XX$ and $\gamma\in \Gamma_{\tilde\eta}$ arbitrary we have by $\Gamma$-equivariance 
$$ \mu_x(\tilde \eta)=\mu_x(\gamma^{-1}\tilde\eta)=\mu_{\gamma x}(\tilde\eta).$$
From the assumption that $\tilde\eta$ is a point mass and 
property (iii) in Definition~\ref{bdensi} we get
$$ 1=\frac{\mu_{\gamma x}(\tilde\eta)}{\mu_{x}(\tilde\eta)}=\e^{\langle b,\bv_{\tilde\eta}(x,\gamma x)\rangle},$$
so  $\langle b,\bv_{\tilde\eta}(x,\gamma x)\rangle =0$ for all $x\in\XX$ and all $\ging_{\tilde\eta}$.
 
Next let $\gamma,\hat \gamma\in\Gamma$ \st
$\gamma\Gamma_{\tilde\eta}=\hat\gamma\Gamma_{\tilde\eta}\in \Gamma/\Gamma_{\tilde\eta}$. Then
$\hat\gamma^{-1}\gamma\in\Gamma_{\tilde\eta}$ and therefore\break  $\langle b,\bv_{\tilde\eta}(\gamma^{-1} x, (\hat\gamma^{-1}\gamma)\gamma^{-1}x)=0$; using the cocycle identity for the Busemann vector we get
\begin{align}
\langle b,\bv_{\tilde\eta}(x,\gamma^{-1} x)\rangle&=\langle b,\bv_{\tilde\eta}(x,\gamma^{-1} x)\rangle +\langle b,\bv_{\tilde\eta}(\gamma^{-1}x,\hat\gamma^{-1}\gamma \gamma
^{-1}x)\rangle\nonumber\\
&=\langle b,\bv_{\tilde\eta}(x,\hat\gamma^{-1}\gamma \gamma
^{-1}x)\rangle = \langle b,\bv_{\tilde\eta}(x,\hat\gamma^{-1}x)\rangle.\tag*{$\scriptstyle\square$}
\end{align}

\begin{lem}\label{sumconv}
If $\tilde\eta\in
\rand_\theta$ is a point mass for a $( b,\theta)$-density $\mu$, then the sum
$$ \sum \e^{ \langle b,\bv_{\tilde\eta}(\xo,\gamma^{-1}\xo)\rangle}$$
taken over a system of coset representatives of $\Gamma/\Gamma_{\tilde\eta}$
 converges. \end{lem}
\prf. \ \   If $\gamma$ and $\hat\gamma$ are representatives of different cosets in
$\Gamma/\Gamma_{\tilde\eta}$, then
$\gamma\tilde\eta\neq \hat\gamma\tilde\eta$  and hence, by $\Gamma$-equivariance, 
the sum
$\sum\mu_{\gamma^{-1}\xo}(\tilde\eta)=\sum\mu_\xo(\gamma\tilde\eta)$
over a system of coset representatives of $\Gamma/\Gamma_{\tilde\eta}$ is bounded above by
$\mu_o(\rand)$. 
By property (iii) in Definition~\ref{bdensi} and the
assumption that
$\tilde\eta$ is a point mass we conclude that the sum
$$ \sum  \e^{\langle b,\bv_{\tilde\eta}(\xo,\gamma^{-1}\xo)\rangle} =
  \sum \frac{\mu_{\gamma^{-1}\xo} (\tilde\eta)}{\mu_{\xo}(\tilde\eta)}  =
\frac1{\mu_\xo(\tilde\eta)} \sum  \mu_{\gamma^{-1}\xo} (\tilde\eta)     $$
over a system of coset representatives of $\Gamma/\Gamma_\eta$
is bounded above by\break 
\hbox{$\displaystyle\frac{\mu_\xo(\rand)}{\mu_\xo(\tilde\eta)}$}.
Since $\mu_\xo$
is a finite measure and $\mu_\xo(\tilde\eta)>0$, the above sum
converges.\qed

\begin{thr}\label{radpoint}
If $\delta_\theta(\Gamma)>0$ then a regular radial limit point $\tilde\eta\in \radlim\cap \regrand$ is not a
point mass for any \bd y. 
\end{thr}
\prf. \ \  Let $\mu$ be a \bd y. If $\tilde\eta\notin
\rand_\theta$, then $\tilde\eta\notin\supp (\mu_\xo)$, hence $\tilde\eta$ cannot be a point
mass. In particular, it suffices to consider $\theta\in E^+$.

Suppose $\tilde\eta=(\eta_1,\eta_2,\ldots,\eta_r,\theta)\in\radlim\cap \rand_\theta\subset\regrand$ is a point mass for $\mu$. Then by Theorem~\ref{posmeasrad}   we have
$\ \langle b,\theta\rangle=\delta_\theta(\Gamma)>0$, hence by continuity of the map $E\to\RR,\ \hat\theta\mapsto \langle b,\hat\theta\rangle\,$ there exists $\eps>0$ \st  every $\hat\theta\in E$ with $\Vert\hat\theta-\theta\Vert<\eps$ satisfies 
\[ \langle b,\hat\theta\rangle \ge \frac{\delta_\theta(\Gamma)}2>0. \]
 Moreover, by the formula~(\ref{radlimthetacover}) for the radial limit set in $\rand_\theta$ there exists  a constant
$c>0$ and a sequence $(\gamma_n)=\big((\gamma_{n,1}, \gamma_{n,2},\ldots,\gamma_{n,r})\big)\subset\Gamma$  \st 
$\Vert\thet(\xo,\gamma_n\xo)-\theta\Vert <\eps$ 
and $\tilde\eta\in S(o: B_{\gamma_n o}(c))$ for all $n\in\NN$.
Corollary \ref{esti}  implies
\[\bs_{\eta_i}(o_i,\gamma_{n,i} o_i)>d_i(o_i,\gamma_{n,i} o_i)-2c\quad\text{for all}\quad  n\in\NN\quad\text{and all}\quad  i\in\{1,2,\ldots,r\},\]
and by choice of $\eps>0$ we have $\langle b,\thet(\xo,\gamma_n\xo)\rangle \ge q$ for all $n\in\NN$. 
Summarizing we  conclude
$$ \langle b,\bv_{\tilde\eta}(\xo,\gamma_n \xo)\rangle>\langle b, H(o,\gamma_n o) \rangle-2 \Vert b\Vert c\ge q\cdot d(\xo,\gamma_n\xo)-2 \Vert b\Vert c \  \longrightarrow\  \infty\quad\text{as}\quad n\to\infty.$$
 Passing to a subsequence if necessary we may therefore
assume that 
$\ \langle b,\bv_{\tilde\eta}(\xo,\gamma_{n} \xo)\rangle\  $
is strictly increasing to infinity as $n\to\infty$.

Now suppose there exist
$l,j \in\NN$, $l\ne j$ \st 
$\gamma_l^{-1}\Gamma_{\tilde\eta} = \gamma_j^{-1}\Gamma_{\tilde\eta}$.
Since $\tilde\eta$ is a point mass for $\mu$,  Lemma~\ref{busnull} implies
$$\langle b,\bv_{\tilde\eta}(\xo,\gamma_{j} \xo)\rangle =\langle b,\bv_{\tilde\eta}(\xo,\gamma_{l} \xo)\rangle,$$
in contradiction to the choice of the subsequence $(\gamma_{n})$. 
Hence
\hbox{$\gamma_l^{-1}\Gamma_{\tilde\eta}{\neq}\gamma_j^{-1}\Gamma_{\tilde
\eta}$} for all $l\ne j$, and
the sum $\displaystyle\sum \e^{\langle b,\bv_{\tilde\eta}(\xo,\gamma\xo)\rangle}$ over a system of coset representatives of $\Gamma/\Gamma_{\tilde\eta}$ is
bounded below by 
$$\sum_{n\in\NN}\e^{\langle b,\bv_{\tilde\eta}(\xo,\gamma_{n} \xo)\rangle} .$$
The divergence of this series yields a contradiction to {\rm Lemma}~\ref{sumconv}  and we conclude
that  $\tilde\eta$ cannot be a point mass for
$\mu_\xo$.\qed

\bibliography{Bibliographie} 

\vspace{1.5cm}
\noindent Gabriele Link\\
Institut f\"ur Algebra und Geometrie\\
Karlsruher Institut f\"ur Technologie (KIT)\\
Kaiserstr. 89-93  \\[1mm]
D-76133 Karlsruhe\\
e-mail:\ gabriele.link@kit.edu

\end{document}